\documentclass{amsart}
\usepackage{amssymb}
\usepackage{amsmath, amscd}
\usepackage{amsthm}
\swapnumbers
\usepackage{mathrsfs}
\usepackage{perpage}
\usepackage{mathdots}
\usepackage[all,cmtip]{xy}
\usepackage{setspace}
\usepackage{amsbsy}
\usepackage{chngcntr}
\usepackage{bm}
\usepackage[foot]{amsaddr}
\usepackage{url}
\usepackage[T1]{fontenc}
\usepackage{verbatim}
\usepackage{mathrsfs}
\usepackage{ifthen}
\usepackage{manfnt}
\usepackage{blkarray}
\usepackage{multirow}
\usepackage{enumitem}
 \usepackage{xr-hyper}
\usepackage{hyperref}
\usepackage{etoolbox}
\usepackage[left=2cm,head=13pt,
top=1.5cm,right=2cm, bottom=1.5cm]{geometry}
\usepackage{tikz-cd}

\setenumerate[0]{label=(\alph*)}

\title{Hodge-Chern classes and strata-effectivity in tautological rings}
\AtEndDocument{\bigskip{\footnotesize%

\noindent
Simon Cooper \textsc{Department of Mathematics, Stockholm University} \par
 \noindent \texttt{cooper.simonjohn@gmail.com}
 
\noindent Wushi Goldring
  \textsc{Department of Mathematics, Stockholm University} \par
 \noindent \texttt{wushijig@gmail.com} \par
 }}

\newtheorem{theorem}[subsubsection]{Theorem}
\newtheorem{lemma}[subsubsection]{Lemma}
\newtheorem{conjecture}[subsubsection]{Conjecture}

\newtheorem{corollary}[subsubsection]{Corollary}
\newtheorem{definition}[subsubsection]{Definition}
\newtheorem{question}[subsubsection]{Question}
\newtheorem{remark}[subsubsection]{Remark}
\newtheorem{proposition}[subsubsection]{Proposition}

\newtheorem{theme}[subsubsection]{Theme}

\theoremstyle{remark}

\newtheorem{rmk}[subsubsection]{Remark}

\newskip\procskipamount
\newskip\interskipamount
\newskip\refskipamount
\newcommand{\gofqp}{\GG(\qp)}

\newcommand{\gofafp}{\GG(\AA_f^p)}

\newcommand{\GZip}{\mathop{\text{$G$-{\tt Zip}}}\nolimits}
\newcommand{\tildeGZip}{\mathop{\text{$\tilde G$-{\tt Zip}}}\nolimits}
\newcommand{\word}{\mathop{\text{$w$-{\tt ord}}}\nolimits}
\newcommand{\HZip}{\mathop{\text{$H$-{\tt Zip}}}\nolimits}

\newcommand{\GoneZip}{\mathop{\text{$G_1$-{\tt Zip}}}\nolimits}
\newcommand{\GtwoZip}{\mathop{\text{$G_2$-{\tt Zip}}}\nolimits}
\newcommand{\GiZip}{\mathop{\text{$G_i$-{\tt Zip}}}\nolimits}

\newcommand{\GdjZip}{\mathop{\text{$G_{d_j}$-{\tt Zip}}}\nolimits}

\newcommand{\GLnZip}{\mathop{\text{$\GL(n)$-{\tt Zip}}}\nolimits}

\newcommand{\GF}{\mathop{\text{$G$-{\tt ZipFlag}}}\nolimits}
\newcommand{\HF}{\mathop{\text{$H$-{\tt ZipFlag}}}\nolimits}

\newcommand{\gal}{\mathsf{Gal}}
\DeclareMathOperator{\im}{im}
\DeclareMathOperator{\Grif}{Grif}
\DeclareMathOperator{\grif}{grif}

\newcommand{\ad}{\mathsf{ad}}

\newcommand{\aut}{\mathsf{Aut}}

\newcommand{\Sbt}{\mathsf{Sbt}}

\newcommand{\pha}{\mathsf{pHa}}
\newcommand{\type}{\mathsf{type}}
\newcommand{\cent}{\textnormal{Cent}}

\newcommand{\Zero}{\textnormal{Zero}}
\newcommand{\inj}{\textnormal{inj}}

\newcommand{\Ind}{\mathsf{Ind}}

\newcommand{\Pic}{\mathsf{Pic}}

\newcommand{\sgn}{\textnormal{sgn}}

\newcommand{\Acal}{{\mathcal A}}
\newcommand{\Bcal}{{\mathcal B}}

\newcommand{\Lcal}{{\mathcal L}}

\newcommand{\Ocal}{{\mathcal O}}
\newcommand{\Pcal}{{\mathcal P}}

\newcommand{\Vcal}{{\mathcal V}}

\newcommand{\Xcal}{{\mathcal X}}
\newcommand{\Ycal}{{\mathcal Y}}
\newcommand{\Zcal}{{\mathcal Z}}

\renewcommand{\AA}{\mathbf{A}}
\newcommand{\BB}{\mathbf{B}}
\newcommand{\CC}{\mathbf{C}}

\newcommand{\FF}{\mathbf{F}}
\newcommand{\GG}{\mathbf{G}}

\newcommand{\NN}{\mathbf{N}}

\newcommand{\PP}{\mathbf{P}}
\newcommand{\QQ}{\mathbf{Q}}
\newcommand{\RR}{\mathbf{R}}

\newcommand{\XX}{\mathbf{X}}

\newcommand{\ZZ}{\mathbf{Z}}

\newcommand{\Cscr}{{\mathscr C}}

\newcommand{\Lscr}{{\mathscr L}}
\newcommand{\Mscr}{{\mathscr M}}

\newcommand{\Sscr}{{\mathscr S}}

\newcommand{\Vscr}{{\mathscr V}}
\newcommand{\Wscr}{{\mathscr W}}

\newcommand{\csf}{{\mathsf c}}

\newcommand{\hsf}{{\mathsf h}}

\newcommand{\Asf}{\mathsf{A}}
\newcommand{\Bsf}{\mathsf{B}}
\newcommand{\Csf}{\mathsf{C}}
\newcommand{\Dsf}{\mathsf{D}}
\newcommand{\Esf}{\mathsf{E}}
\newcommand{\Fsf}{\mathsf{F}}
\newcommand{\Gsf}{\mathsf{G}}
\newcommand{\Hsf}{\mathsf{H}}

\newcommand{\Tsf}{\mathsf{T}}

\newcommand{\Xsf}{\mathsf{X}}

\newcommand{\Zsf}{\mathsf{Z}}

\newcommand{\po}{\PP^1}
\newcommand{\poc}{\po(\CC)}

\newcommand{\fp}{\FF_p}
\newcommand{\fpd}{\FF_{p^{d}}}

\newcommand{\leftexp}[2]{{\vphantom{#2}}^{#1}{#2}}

\newcommand{\CharT}{\mathsf{X}^{*}(T)}

\newcommand{\qgeqz}{\QQ_{\geq 0}}

\newcommand{\qp}{\QQ_p}

\newcommand{\std}{{\rm Std}}

\newcommand{\gx}{(\GG, \XX)}

\newcommand{\res}{{\rm Res}}

\newcommand{\der}{{\rm der}}

\newcommand{\End}{{\rm End}}

\newcommand{\af}{\AA_f}

\newcommand{\id}{{\rm Id}}

\newcommand{\Th}{{\rm Th.}}

\newcommand{\Rmk}{{\rm Rmk.}}

\newcommand{\Lem}{{\rm Lem.}}

\newcommand{\Prop}{{\rm Prop.}}

\newcommand{\Exa}{{\rm Exa.}}

\newcommand{\cf}{{\em cf. }}
\newcommand{\ie}{i.e.,\ }

\newcommand{\eg}{e.g.,\ }

\newcommand{\fil}{{\rm Fil}}
\newcommand{\Gr}{{\rm Gr}}
\newcommand{\gr}{{\rm gr}}

\newcommand{\can}{{\rm can}}

\newcommand{\Eff}{\mathsf{Eff}}
\newcommand{\strat}{\mathsf{strat}}
\newcommand{\taut}{\mathsf{taut}}

\newcommand{\dR}{{\rm dR}}

\newcommand{\CH}{\Csf\Hsf}
\newcommand{\CHQ}{\Csf\Hsf_{\QQ}}
\newcommand{\CHQs}{\Csf\Hsf^*_{\QQ}}

\newcommand{\iw}{\leftexp{I}{W}}

\renewcommand{\div}{{\rm div}}

\newcommand{\GL}{\mathsf{GL}}
\newcommand{\PGL}{\mathsf{PGL}}
\newcommand{\SL}{\mathsf{SL}}

\newcommand{\Sp}{\mathsf{Sp}}
\newcommand{\GSp}{\mathsf{GSp}}
\newcommand{\Spin}{\mathsf{Spin}}

\newcommand{\U}{\mathsf{U}}

\author{Simon Cooper and Wushi Goldring}
\date{\today}
\makeatletter
\let\c@equation=\c@subsubsection

\let\c@figure=\c@subsubsection

\begin{document}
\pagestyle{plain}

\begin{abstract}
Given a connected, reductive $\fp$-group $G$, a cocharacter $\mu \in \Xsf_*(G)$ and a smooth zip period map $\zeta:X \to \GZip^{\mu}$, we study which classes in the Wedhorn-Ziegler tautological rings $\Tsf^*(X), \Tsf^*(Y)$ of $X$ and its flag space $Y \to \GF^{\mu}$ are \textit{strata-effective}, meaning that they are non-negative rational linear combinations of pullbacks of classes of zip (flag) strata closures. 
Two special cases are: (1) When $X=\GZip^{\mu}$ and the tautological rings $\Tsf^*(X)=\CHQ(\GZip^{\mu})$, $\Tsf^*(Y)=\CHQ(\GF^{\mu})$ are the entire Chow ring, and (2) When $X$ is the special fiber of an integral canonical model of a Hodge-type Shimura variety -- in this case the strata are also known as Ekedahl-Oort strata. 
We focus on the strata-effectivity of three types of classes: (a) Effective tautological classes, (b) Chern classes of Griffiths-Hodge bundles and (c) Generically $w$-ordinary curves. We connect the question of strata-effectivity in (a) to the global section `Cone Conjecture' of Goldring-Koskivirta. For every representation $r$ of $G$, we conjecture that the Chern classes of the Griffiths-Hodge bundle associated to $(G, \mu,r)$ are all strata-effective. This provides a vast generalization of a result of Ekedahl-van der Geer that the Chern classes of the Hodge vector bundle on the moduli space of principally polarized abelian varieties $\Acal_{g,\FF_p}$ in characteristic $p$ are represented by the closures of $p$-rank strata. We prove several instances of our conjecture, including the case of Hilbert modular varieties, where the conjecture says that all monomials in the first Chern classes of the factors of the Hodge vector bundle are strata-effective. We prove results about each of (a), (b), (c) which have applications to Shimura varieties and also in cases where no Shimura variety exists.  
\end{abstract}
\maketitle
\tableofcontents
\section{Introduction}
\label{sec-intro}
The aim of this paper is to unify and generalize three seemingly different threads (\S\ref{sec-3-directions}) concerning moduli spaces $X$ in characteristic $p>0$ by viewing them within the frameworks of \textit{$G$-Zip geometricity} (\S\ref{sec-zip-geom}) and \textit{generation by groups}. The former was laid out by one of us (W.G.) with J.-S. Koskivirta in \cite{Goldring-Koskivirta-Strata-Hasse} and developed further in \cite{Goldring-Koskivirta-global-sections-compositio}; for the latter  see \cite{Goldring-griffiths-bundle-group-theoretic}. To this end, the main conceptual idea of this paper is that several different kinds of questions about such moduli spaces fit into the setting of \textit{strata-effectivity in tautological rings} (\S\ref{sec-def-strata-effective}). An important special case is when $X$ arises from a Hodge-type Shimura variety (\S\ref{sec-intro-Hodge-type}). However, guided by these frameworks, one of the paper's main goals is to lay out a theory, questions and conjectures about strata-effective classes in a far broader context than Shimura varieties.
\subsection{\texorpdfstring{$G$}{G}-Zip Geometricity}
\label{sec-zip-geom}
Let $p$ be a prime and let $k$ be an algebraic closure of $\fp$. Let $G$ be a connected, reductive $\fp$-group and let $\mu \in \Xsf_*(G)$ be a cocharacter over $k$. Associated to the pair $(G,\mu)$, Pink-Wedhorn-Ziegler \cite{Pink-Wedhorn-Ziegler-zip-data,Pink-Wedhorn-Ziegler-F-Zips-additional-structure} define a stack $\GZip^{\mu}$. Morphisms $X \to \GZip^{\mu}$ give a mod $p$ analogue of Hodge structures with $G$-structure. For example, given a sufficiently nice family of $k$-schemes $Y/X$ and a non-negative integer $i$, the de Rham cohomology $H^i_{\dR}(Y/X)$ is classified by a morphism $X \to \GLnZip^{\mu}$, where $n$ is the rank of $H^i_{\dR}(Y/X)$ and $\mu$ is deduced from the Hodge filtration, just like the Hodge-Deligne cocharacter in classical Hodge theory. For more details on this analogy, \cf the introduction of \cite{Goldring-Koskivirta-global-sections-compositio}. Motivated by the analogy with classical Hodge theory, we refer to morphisms of $k$-stacks $\zeta:X \to \GZip^{\mu}$ as \textit{Zip period maps} (including the case $X=\GZip^{\mu}$, $\zeta=\id$).  

Recall from \cite[Question~A]{Goldring-Koskivirta-global-sections-compositio} that the question driving $G$-Zip geometricity is:
\begin{question}
\label{q-G-Zip-geom}
    Given a Zip period map 
\addtocounter{equation}{-1}
\begin{subequations}
\begin{equation}
 \label{eq-intro-zip-period}   
\zeta: X \to \GZip^{\mu},
\end{equation}    
\end{subequations}
 what geometry of $X$ is determined by properties of the morphism $\zeta$ and the stack $\GZip^{\mu}?$
\end{question}
\subsubsection{Examples of \texorpdfstring{$G$}{G}-Zip Geometricity}
Some works illustrating
the philosophy of $G$-Zip geometricity are:
\begin{enumerate}
\item The original paper of Moonen-Wedhorn \cite{Moonen-Wedhorn-Discrete-Invariants} on $F$-Zips,
\item
\label{item-past-strata-hasse}
The works of Goldring-Koskivirta on strata Hasse invariants and their applications to the Langlands correspondence \cite{Goldring-Koskivirta-Strata-Hasse} and on the cone of global sections \cite{Goldring-Koskivirta-global-sections-compositio}, \cite{Goldring-Koskivirta-GS-cone}, \cite{Goldring-Koskivirta-divisibility};
\item The work of Brunebarbe-Goldring-Koskivirta-Stroh on ampleness of automorphic bundles \cite{Brunebarbe-Goldring-Koskivirta-Stroh-ampleness};
\item The works of Wedhorn-Ziegler \cite{Wedhorn-Ziegler-tautological} and one of us (S. C.) \cite{Cooper-tautological-ring-hilbert} on tautological subrings of Chow rings;
\item
\label{item-past-Ogus}
The work of Reppen \cite{Reppen-ogus-hilbert} and Goldring-Reppen on a  generalization of Ogus' Principle relating the Hasse invariant's vanishing order and `Frobenius and the Hodge filtration' to Zip period maps \cite{Goldring-Reppen}.
\end{enumerate}
\subsubsection{Reduction to flag varieties}
\label{sec-reduction-flag-varieties}
Following the analogy with classical Hodge theory, the works above as well as this paper also develop and exploit a second reduction from $\GZip^{\mu}$ to flag varieties. In~\ref{item-past-strata-hasse}-\ref{item-past-Ogus} as well as here, the reduction to flag varieties is based on the following diagram introduced in \cite{Goldring-Koskivirta-Strata-Hasse}, whereby a flag space $\GF^{\mu}$ over $\GZip^{\mu}$ is used to connect $\GZip^{\mu}$ to the Schubert stack $\Sbt$, for a well-chosen $\fp$-Borel $B$ of $G$:
$$\xymatrixrowsep{1pc} \xymatrix{ & \GF^{\mu} \ar[rd]^-{\hsf} \ar[ld] & \\ \GZip^{\mu} & & \Sbt:=[B \backslash G/B],  }
$$
\subsubsection{A key test-case: Hodge-type Shimura varieties} 
\label{sec-intro-Hodge-type}
Let $\gx$ be a Hodge-type Shimura datum. Assume that $\GG$ is unramified at $p$. Let $K_p \subset \gofqp$ be a hyperspecial maximal compact subgroup. By the work of Kisin \cite{Kisin-Hodge-Type-Shimura} and Vasiu \cite{Vasiu-Preabelian-integral-canonical-models}, as $K^p$ ranges over open, compact subgroups of $\gofafp$,  the associated projective system of Shimura varieties  admits an integral canonical model $(\Sscr_{K_pK^p}\gx)_{K^p}$ in the sense of Milne \cite{Milne-integral-canonical-models}. Set $K:=K_pK^p$ and let $S_K$ be the special $k$-fiber of $\Sscr_{K_pK^p}\gx$. Given a symplectic embedding of $\gx$ into a Siegel-type datum $(\GSp(2g), \XX_g)$, for all sufficiently small $K^p$ there exists a level $K' \subset \GSp(2g, \af)$ and an induced finite map from $S_K$ to the special $k$-fiber of the Siegel-type Shimura variety $S_{g,K'}$ (\cf \cite[(2.3.3)]{Kisin-Hodge-Type-Shimura}, \cite[\S4.1]{Goldring-Koskivirta-Strata-Hasse} for more details). If $Y/S_K$ is the resulting family of abelian schemes, the Zip period map associated to $H^1_{\dR}(Y/S_K)$ factors through a smooth (Zhang \cite{Zhang-EO-Hodge}) surjective (Kisin-Madapusi Pera-Shin \cite{Kisin-Madapusi-Pera-Shin-Honda-Tate}) morphism 
\addtocounter{equation}{-1}
\begin{subequations}
\begin{equation}
\label{eq-zeta-shimura}
\zeta:S_K \to \GZip^{\mu},
\end{equation}    
\end{subequations}
 where $G$ is the reductive $\fp$-group deduced from the $\QQ$-group $\GG$ and $\mu \in \Xsf_*(G)$ is a representative of the conjugacy class of cocharacters deduced from the Hermitian symmetric space $\XX$. 
Associated to a suitable choice of combinatorial data $\Sigma$, let $S_K^{\Sigma}$ be the smooth toroidal compactification of $S_K$  constructed by Madapusi Pera \cite{Madapusi-Hodge-tor}. For sufficiently nice $\Sigma$, $\zeta$ extends to $\zeta^{\Sigma}:S_K^{\Sigma} \to \GZip^{\mu}$ by Goldring-Koskivirta \cite[\S6]{Goldring-Koskivirta-Strata-Hasse} and $\zeta^{\Sigma}$ is smooth by Andreatta \cite{Andreatta-modp-period-maps} (previously shown by Lan-Stroh \cite{Lan-Stroh-stratifications-compactifications} in the PEL case).
\subsection{Three threads to unify and generalize}
\label{sec-3-directions}
Consider the three threads~\ref{sec-intro-ekedahl-geer}-\ref{sec-intro-wedhorn-ziegler}: 
\subsubsection{Chern classes of the Hodge vector bundle on \texorpdfstring{$S_{g,K'}$}{!} are positive multiples of $p$-rank strata classes, d'apr\`es Ekedahl-van der Geer \cite{Ekedahl-Geer-EO}}
\label{sec-intro-ekedahl-geer}
In the Siegel case $S_K=S_{g,K'}$, \cite[\Th~12.4]{Ekedahl-Geer-EO} shows  that the Chern classes of the 
Hodge vector bundle $\Omega^{\can}:=\fil^1H^{1,\can}_{\dR}$ are all \textit{effective} in the Chow ring with rational coefficients $\CHQ^*(S_{g,K'}^{\Sigma})$, where $H^{1,\can}_{\dR}(S_{g, K'})$ denotes the canonical extension of the de Rham bundle $H^{1}_{\dR}(Y/S_{g, K'})$ to $S_{g,K'}^{\Sigma}$. More precisely: For every $i$, the $i$th Chern class $\csf_i(\Omega^{\can})$ is a strictly positive rational multiple of the cycle class of the locus of abelian varieties of $p$-rank $\leq g-i$.   
\subsubsection{The global section `Cone Conjecture' on \texorpdfstring{$G$}{G}-Zip schemes, d'apr\texorpdfstring{\`}{!}es Goldring-Koskivirta \texorpdfstring{\cite{Goldring-Koskivirta-global-sections-compositio,Goldring-Koskivirta-GS-cone,Goldring-Koskivirta-divisibility,Goldring-Koskivirta-rank-2-cones}}{!}}
\label{sec-intro-global-sections}
\addtocounter{equation}{-1}
\begin{subequations}
Let $P \subset G_k$ be the parabolic subgroup of non-positive $\mu$-weights. As in \cite[1.2.3]{Goldring-Koskivirta-Strata-Hasse}, assume for simplicity that there exists an $\fp$-Borel $B$ of $G$ contained in $P$. 
Since $\GZip^{\mu}=[E\backslash G]$ where $E$ is the zip group of $(G,\mu)$,  every $k$-representation $r$ of $P$ admits an associated `automorphic' vector bundle $\Wscr(r)$ on $\GZip^{\mu}$ via the projection $E \to P$. 
Given a character $\lambda \in \Xsf^*(B)$, let $\Vscr(\lambda)$ be the automorphic vector bundle associated to the induced (Borel-Weil) representation $\Ind_{B}^P \lambda=H^0(P/B, \Lscr(\lambda))$, where $\Lscr(\lambda)$ is the associated line bundle on the flag variety $P/B$. 
For a Hodge-type Shimura variety~\eqref{eq-zeta-shimura}, the pullbacks $\zeta^*\Vscr(\lambda)$ are precisely reductions mod $p$ of canonical integral models of the classical automorphic vector bundles on the Shimura variety over $\CC$, whose coherent cohomology gives automorphic forms. 

Given a zip period map~\eqref{eq-intro-zip-period}, 
the \textit{global section cone} of $X$ is 
\begin{equation}
    \Cscr_X=\{\lambda \in \Xsf^*(B) \ | \ H^0(X, \zeta^*\Vscr(n\lambda)) \neq 0 \textnormal{ for some } n>0\}.
\end{equation}
The `Cone Conjecture', originally stated in \cite{Goldring-Koskivirta-global-sections-compositio} and subsequently generalized in \cite{Goldring-Koskivirta-rank-2-cones} and \cite{Goldring-Koskivirta-GS-cone} states that, under mild assumptions on $\zeta$ and $\mu$ (see \ref{conj-cone}), the global section cones of $X$ and $\GZip^{\mu}$ are equal: 
\begin{equation}
\label{eq-intro-cone-conj}    
\Cscr_X=\Cscr_{\GZip^{\mu}}.
\end{equation}
 In terms of the basic $G$-Zip Geometricity Question~\ref{q-G-Zip-geom}, the Cone Conjecture says that the global section cone $\Cscr_X$ is a geometric invariant of $X$ which is determined by the group theory of $(G,\mu)$ and the regularity of $\zeta$. 
\end{subequations}
\subsubsection{Tautological rings of mod \texorpdfstring{$p$}{!} Shimura varieties via \texorpdfstring{$G$}{G}-Zips, d'apr\texorpdfstring{\`}{!}es Wedhorn-Ziegler 
\cite{Wedhorn-Ziegler-tautological}}
\label{sec-intro-wedhorn-ziegler}
\addtocounter{equation}{-1}
\begin{subequations}
Given a smooth Zip period map~\eqref{eq-intro-zip-period}, Wedhorn-Ziegler define the \textit{tautological ring} $\Tsf^*(X) \subset \CHQ^*(X)$ to be the subring generated by all Chern classes of all automorphic bundles $\zeta^*\Vscr(\lambda)$; equivalently $\Tsf^*(X)=\zeta^*\CHQ^*(\GZip^{\mu})$ is the  image of the pullback along $\zeta$. When $X=S_K^{\Sigma}$ is a toroidal compactification of a Hodge-type Shimura variety~\eqref{sec-intro-Hodge-type}, they show that $\zeta^*$ is injective; hence 
\begin{equation}
\label{eq-Wedhorn-Ziegler-iso}
\zeta^*:\CHQ^*(\GZip^{\mu}) \stackrel{\sim}{\to} \Tsf^*(S_K^{\Sigma}).    
\end{equation} 
In terms of our guiding Question~\ref{q-G-Zip-geom}, the Wedhorn-Ziegler isomorphism~\eqref{eq-Wedhorn-Ziegler-iso} says that, when $X=S_K^{\Sigma}$, the answer includes the tautological ring. 
Highlighting the theme~\ref{sec-reduction-flag-varieties} of successive reduction from $X$ to $\GZip^{\mu}$ to the flag variety compact dual $G/P$, Brokemper and Wedhorn-Ziegler also show that $\CHQs(\GZip^{\mu}) \cong H^{2*}(\GG_{\CC}/\PP_{\CC}, \QQ)$, where $\GG_{\CC}/\PP_{\CC}$ is the flag variety of the same type~(\ref{sec-chow-zip-flag}-\ref{cor-chow-zip-flag}).
\end{subequations}
\subsection{Strata-effective classes}
\label{sec-def-strata-effective}
Let $X$ be a stack endowed with a stratification
\begin{equation}
\label{eq-intro-strat}
    X=\bigsqcup_{w \in I}X_w
\end{equation} by locally closed substacks $X_w$ indexed by a finite set $I$. 
We say that a class $\eta \in \CHQ^*(X)$ is \underline{strata-effective} if 
\begin{equation}
\label{eq-def-strata-effective}
\eta=\sum_{w \in I} a_w [\overline{X}_w] \quad \textnormal{ with } a_w \in \qgeqz  
\end{equation}
is a non-negative $\QQ$-linear combination of classes of strata closures. In particular, a strata-effective class is effective.

\subsubsection{Flag spaces}
\label{sec-intro-flag-space}
To illustrate the broader scope of strata-effectivity in the context of Zip period maps, recall (\cite{Goldring-Koskivirta-Strata-Hasse} \cite{Goldring-Koskivirta-zip-flags}) that $\GZip^{\mu}$ admits a tower of partial flag spaces $\GF^{\mu} \to \GF^{\mu,P_0} \to \GZip^{\mu}$ indexed by intermediate parabolics $B \subset P_0 \subset P$. 
The quotient stacks $\GF^{\mu}$,  $\GF^{\mu, P_0}$ and $\GZip^{\mu}$ all admit stratifications~\ref{sec-fine-strat}. 
Let $\zeta_{P_0}:X^{P_0} \to \GF^{\mu,P_0}$ and $\zeta_Y: Y \to \GF^{\mu}$ be the base changes of $\zeta$ along $\GF^{\mu,P_0} \to \GZip^{\mu}$ and $\GF^{\mu} \to \GZip^{\mu}$. If $\zeta$ is smooth, then the stratifications pull back to ones of $X$ and $Y$. When $X$ is the special fiber of a Hodge-type Shimura variety as in~\ref{sec-intro-Hodge-type}, the resulting stratification of $X$ is often called the Ekedahl-Oort (EO) stratification.  Brokemper \cite{Brokemper-tautological} and Wedhorn-Ziegler \cite{Wedhorn-Ziegler-tautological} show that the classes of strata-closures give  bases of the vector spaces $\CHQ^*(\GZip^{\mu})$ and $\CHQ^*(\GF^{\mu})$; if $\zeta$ is smooth, the same holds for the tautological rings $\Tsf^*(X)$ and $\Tsf^*(Y)$. 
\begin{theme}
This paper studies strata-effectivity in the context of smooth Zip period maps~\eqref{eq-intro-zip-period} and their base changes to (partial) flag spaces. Strata-effectivity is used to unify and generalize the three threads highlighted in~\ref{sec-3-directions}.    
\end{theme}  
A basic question about strata-effectivity is:
\begin{question}[Effective tautological vs tautologically effective classes]
\label{q-eff-taut-taut-eff}
Suppose $\eta$ is a tautological class in  $\Tsf^i(X)$ or $\Tsf^i(Y)$ for some $i$. If $\eta$ is effective (as a class in Chow), is $\eta$ also strata-effective?  
\end{question}
Since the classes of strata-closures are `tautologically effective',~\ref{q-eff-taut-taut-eff} is an attempt to compare effective, tautological classes with `tautologically effective' ones. Note that~\ref{q-eff-taut-taut-eff} is already highly nontrivial for $\Tsf^i(\GF^{\mu})=\CHQ^i(\GF^{\mu})$ when $\zeta$ and $\zeta_Y$ are the identity.
\subsubsection{Example: The Schubert-Bruhat-Chevalley stratification of a flag variety \texorpdfstring{$G/B$}{!}}
\label{ex-flag-variety-G/B}
The best-case scenario concerning~\ref{q-eff-taut-taut-eff} is given by the example where $X:=G/B$ is the flag variety of $G$, $I=W$ is the Weyl group of $G$ relative some maximal torus $T \subset B \subset G$ and $X_w:=B\dot wB/B$ are the \textit{Schubert cells}
of $X$. Since Schubert cells are affine spaces, the classes of Schubert varieties $([\overline{X}_w])_{w \in W}$ form a basis of $\CHQ^*(X)=H^*(\GG_{\CC}/\BB_{\CC}, \QQ)$ \cite[1.9.1,~1.10.2]{Fulton-Intersection-theory-book}. It is an easy consequence of the existence of opposite Schubert cells that every effective class in $\CHQ^*(X)$ is strata-effective \cite[1.3.6]{brion-lectures-flag-varieties}.     

\subsubsection{\texorpdfstring{$\textdbend$}{!} Effective classes need not be strata-effective}
One of the basic subtleties regarding the `strata-effective' notion is that a tautological class $\eta \in \Tsf^*(X)$ which is effective as an element of $\CHQ^*(X)$ need not be strata-effective: One may have $\eta=[Z]$ for some subscheme $Z \subset X$ and at the same time a decomposition as in~\eqref{eq-def-strata-effective} but with some negative coefficients $a_w<0$. For an example, see \ref{sec-efftaut-notstrataeff}.

On the other hand, there are interesting cases where this does hold. We use a connection between the Cone Conjecture and strata-effectivity to illustrate many examples of both effective, tautological classes which are strata-effective and ones which are not. It would be interesting to consider a more sophisticated version of~\ref{q-eff-taut-taut-eff} where the collection of `tautologically effective' classes is enlarged to include other classes beyond strata-closures which arise naturally on $\GZip^{\mu}$ and $\GF^{\mu}$.    

\subsection{Strata-effectivity of Hodge-Chern classes}
\label{sec-intro-Hodge-chern}
Since the $p$-rank is constant on an Ekedahl-Oort stratum, the result of Ekedahl-van der Geer recalled in~\ref{sec-intro-ekedahl-geer} implies the following about strata-effectivity: 
The Chern classes of the Hodge vector bundle of the Siegel variety $S_{g,K'}$ are all strata-effective. 
We propose a conjecture generalizing this statement in several different directions: 
Let $r:G \to \GL(V)$ be a $k$-representation of $G$. Dropping the assumptions in \cite[\S\S3.1-3.2]{Goldring-griffiths-bundle-group-theoretic} that $r$ is an $\fp$-representation gives the\textit{ Griffiths-Hodge vector bundle} $\Grif(\GZip^{\mu},r)$ on $\GZip^{\mu}$ (\ref{sec-griffiths-bundle}). Pullback along a Zip period map~\eqref{eq-intro-zip-period} gives the Griffiths-Hodge vector bundle $\Grif(X,\zeta,r)$ on $X$. Highlighting the analogy between $G$-Zips and classical Hodge theory, the bundles $\Grif(\GZip^{\mu},r)$ and $\Grif(X,\zeta,r)$ are the Zip period map analogues of the ones studied by Griffiths \cite{Griffiths-IHES-period-integrals} on a Griffiths-Schmid manifold $\Gamma \backslash D$ and on the base $S$ of a $\QQ$-variation of Hodge structure (VHS) by pullback along a classical period map $\Phi: S \to \Gamma \backslash D$. Using these bundles, Griffiths showed that the image $\Phi(S)$ is projective when $S$ is.   
\begin{conjecture}
\label{conj-griffiths-strata-eff} For every triple $(G,\mu,r)$, all the Chern classes of $\Grif(\GZip^{\mu},r)$ are strata-effective.
\end{conjecture}
If $\zeta$ is smooth and surjective, then~\ref{conj-griffiths-strata-eff} implies that all Chern classes of $\Grif(X,\zeta, r)$ are also strata-effective. We call the Chern classes of $\Grif(\GZip^{\mu},r)$ and $\Grif(X,\zeta,r)$ \textit{Hodge-Chern classes}. In the more classical setting of the Griffiths-Hodge vector bundle of a $\QQ$-Variation of Hodge structure over a quasi-projective variety $S$ over $\CC$, one obtains Hodge-Chern classes in $\CHQ^*(S)$ in the same way.

The Hodge vector bundle of a Hodge-type Shimura variety $S_K$ or a smooth toroidal compactification $S_K^{\Sigma}$ (\ref{sec-intro-Hodge-type})  relative a symplectic embedding $\psi:\gx \hookrightarrow (\GSp(W), \XX_g)$ is the special case of the Griffiths-Hodge vector bundle where $r=(\std \circ \psi)_1^{\vee}$ is the dual of the special fiber of an integral model of the composition $\GG \stackrel{\psi}{\to} \GSp(W) \stackrel{\std}{\to} \GL(W)$. 
The vector bundle on $S_K$ associated to $r$ is $H^1_{\dR}(Y/S_K)$~\eqref{sec-intro-Hodge-type}.
More intrinsically in terms of $(G,\mu)$, the Hodge vector bundle is the special case of the Griffiths bundle where the set of $\mu$-weights of $r$ is $\{0,1\}$. 

Observe that, even in the special case of the Hodge vector bundle $\psi^*\Omega^{\can}$ on $S_K^{\Sigma}$,~\ref{conj-griffiths-strata-eff} implies more than the strata-effectivity of the Chern classes of $\psi^*\Omega^{\can}$. 
Roughly, it also implies the strata-effectivity of the `natural' factors of the Hodge vector bundle.   
The representation $\std \circ \psi_1^{\vee}$ of $G$ may be reducible in a non-trivial way over $k$. 
If $r'$ is a sub-quotient of the base-change $(\std \circ \psi)^{\vee}_{1,k}$ to $k$, then~\ref{conj-griffiths-strata-eff} also applies to $(G,\mu,r')$. 
In particular~\ref{conj-griffiths-strata-eff} implies:
\begin{conjecture}
\label{conj-hodge-bundle-strata-eff}    
For every symplectic embedding $\psi$ of Hodge-type Shimura varieties, if $r'$ is a $G$-subrepresentation of $(\std \circ \psi)_1^{\vee}$ with associated bundle $\Wscr(r')$ on $S_K$, then the Chern classes of $\Wscr(r') \cap \psi^*\Omega^{\can}$ are all strata-effective. \end{conjecture}
The strata-effectivity of the Chern classes of the Hodge vector bundle $\psi^*\Omega^{\can}$ does not follow by pullback from the Siegel case: A finite morphism $S_K \to S_{g,K'}$ 
induced by $\psi$
will almost never be flat, so there is no general reason for the pullback of an effective class to be effective. Related to this, for many $\psi$, the preimage of many $p$-rank strata is empty. 

When $S_K^{\Sigma}$ is the $d$-dimensional Hilbert modular variety associated to a degree $d$ totally real extension $F/\QQ$ and $\psi$ is the standard PEL embedding, the Hodge vector bundle $\psi^*\Omega^{\can}=\bigoplus_{i \in \ZZ/d} \omega_i$ decomposes as a direct sum of $d$ line bundles over $k$. Then~\ref{conj-hodge-bundle-strata-eff} predicts that, for every non-empty subset $I \subset \ZZ/d$, the Hodge-Chern monomial $L_I:=\prod_{i \in \ZZ/d}\csf_1(\omega_i)$ is strata-effective. We prove that this prediction is correct. The key is:

\begin{theorem}
\label{th-intro-inert-hilbert}
Let $d \geq 1$ and $G:=\res_{\fp^d/\fp}\GL(2)$. Assume that the parabolic $P$ of negative $\mu$-weights is a Borel subgroup of $G$ and that that $r$ is $k$-subrepresentation of the natural, faithful $2d$-dimensional representation of $G$. 
Then~\ref{conj-griffiths-strata-eff} holds for $(G,\mu, r)$ and there is an explicit formula for the coefficients $a_w$~\eqref{eq-def-strata-effective}. 
\end{theorem}
The explicit formula for the coefficients gives an explicit characterization of the strata-effective classes in the basis of Hodge-Chern monomials of $\Tsf^{|I|}(X)$.
In fact, the coefficients $a_w$ satisfy a reciprocity law, see~\ref{sec-subsubsec-reciprocity}-\ref{prop-reciprocity}.
It is easy to see (\S\ref{sec-functoriality}) that~\ref{conj-griffiths-strata-eff} is stable under products $(G_1,\mu_1)\times (G_2,\mu_2)$ and changing the center of $G$. Hence:
\begin{corollary}
 \label{cor-intro-hilbert}
 Assume that $X=S_K$ is a Hilbert modular variety. Then~\ref{conj-hodge-bundle-strata-eff} holds.
\end{corollary}
Note that, when $G^{\ad}$ is $\fp$-simple as in~\ref{th-intro-inert-hilbert} (equivalently $p$ is inert in $F$), every non-ordinary point of $S_K^{\Sigma}$ has $p$-rank zero. So the effectivity implied by~\ref{cor-intro-hilbert} is nontrivial even for the "full" Chern classes $$\csf_i(\psi^*\Omega^{\can})=\sum_{I \subset \ZZ/d, |I|=i}L_I.$$

As further evidence for~\ref{conj-griffiths-strata-eff} which goes beyond the case of Shimura varieties, building on \cite{Goldring-griffiths-bundle-group-theoretic} and \cite{Goldring-Koskivirta-Strata-Hasse} we show for a general group $G$ that, under a weak `$p$-smallness' assumption introduced in \cite{Goldring-Koskivirta-Strata-Hasse}, all positive powers of the first Chern class of the Griffiths bundle are always strata-effective.
\begin{theorem}[\ref{cor-grif-line-strata-eff}]
\label{th-intro-first-chern-griffiths-effective}
Assume that $r:G \to \GL(V)$ is an $\fp$-representation with central kernel. If $\mu$ is orbitally $p$-close~\ref{def-orb-p-close},  then $\csf_1(\Grif(\GZip^{\mu},r))^m$ is strata-effective for all $m \geq 0$.
\end{theorem}
Let $\tilde{G}$ be the simply-connected cover of the derived subgroup $G^{\der}$ of $G$~\eqref{sec-derived-simply-connected}. 
When the pair $(\type(G), \type(L))$ arises from an $\RR$-group of real rank one, it is shown~(\ref{sec-proportionality}) that all Chern classes of the Griffiths bundle $\Grif(\GZip^{\mu},r)$ are non-negative scalar multiples of the first Chern class $\csf_1(\Grif(\GZip^{\mu}, r)$. Hence: 
\begin{theorem}
\label{th-grif-eff}
Assume that $(\type(G), \type(L))$ is either $(\Xsf_n, \Xsf_{n-1})$ for some $\Xsf \in \{\Asf, \Bsf, \Csf, \Dsf\}$ or $(\Gsf_2, \Asf_{1})$. If $\type(G)=\Gsf_2$, assume $p \geq 5$.
Let
$r$ be a nontrivial, irreducible representation of $\tilde G$ of minimal dimension. Then all the Chern classes of $\Grif(\GZip^{\mu},r)$ are strata-effective.     
\end{theorem}
The possibilities for $r$ are recalled in~\ref{rmk-lowest-fundamental}; when $\type(G) \neq \Asf_n$, ($n \geq 2)$ and $\Dsf_4$, $r$ is the unique smallest fundamental representation.   Pulling back~\ref{th-grif-eff} along the Shimura Zip period map~\eqref{eq-zeta-shimura} shows that all $r$-Hodge-Chern classes are effective for the Shimura  varieties $S_K$ associated to unitary (similitude) groups of signature\footnote{For unitary (similitude) groups of signature $(2,2)$, note the coincidental isomorphism $(\Asf_3, \Asf_1 \times \Asf_1)=(\Dsf_3,\Dsf_2)$.} $(n,1)$, $(2,2)$ and orthogonal/spin (similitude) groups of signature $(2n-2,2)$, $(2n-1,2)$. In the Hodge-type case, for unitary and spin similitude groups of these signatures the same applies to the toroidal compactifications $S_K^{\Sigma}$. By contrast, the cases $(\Csf_n, \Csf_{n-1})$ and $(\Gsf_2, \Asf_1)$ of~\ref{th-grif-eff} lie beyond Shimura varieties.
\subsection{Strata-effectivity \texorpdfstring{\&}{and} the Cone Conjecture}
In the works \cite{Goldring-Koskivirta-global-sections-compositio,Goldring-Koskivirta-rank-2-cones,Goldring-Koskivirta-divisibility,Goldring-Koskivirta-GS-cone,Imai-Koskivirta-zip-schubert} on the Cone Conjecture 
~\eqref{eq-intro-cone-conj}, there is a third cone which plays an important role: The cone $\Cscr_{\pha}$ of \textit{partial Hasse invariants}\footnote{In several of the earlier papers \cite{Goldring-Koskivirta-global-sections-compositio} by Koskivirta and one of us (W.G.), the cone $\Cscr_{\pha}$ was called the Schubert cone and denoted $C_{\Sbt}$.}. 
In the very special case when $X=S_K$ is the special fiber of a Hilbert modular variety~\eqref{sec-intro-Hodge-type}, $\Cscr_{\pha}$ is the cone spanned by the weights of Goren's partial Hasse invariants \cite{Goren-partial-hasse}. 
By construction, one always has the cone inclusions
\begin{equation}
    \Cscr_{\pha} \subset \Cscr_{\GZip^{\mu}} \subset \Cscr_X.
\end{equation}
\begin{theorem}
\label{th-pha-strata-eff}
 For fixed $(G,\mu)$, the following are equivalent:
    \begin{enumerate}
        \item For all $\zeta : X \rightarrow \GZip^\mu$ satisfying~\ref{conj-cone}\ref{item-smooth}-\ref{item-pseudo-complete}, $\Cscr_{X} = \Cscr_{\pha}$.
        \item For all $\zeta : X \rightarrow \GZip^\mu$ satisfying~\ref{conj-cone}\ref{item-smooth}-\ref{item-pseudo-complete}, every effective class $\eta \in \Tsf^1(Y)$ on the flag space $Y/X$~\eqref{sec-intro-flag-space} is strata-effective. 
    \end{enumerate} 
\end{theorem}
In view of~\ref{th-pha-strata-eff}, our basic strata-effectivity Question~\ref{q-eff-taut-taut-eff} can be seen as a generalization of the Cone Conjecture to tautological classes of higher codimension.
In particular, if there is some $\zeta : X \rightarrow \GZip^{\mu}$ for which every effective class $\eta \in \Tsf^1(Y)$ is strata-effective, then the Cone Conjecture holds. Conversely, applying~\ref{th-pha-strata-eff} to cases where $\Cscr_{\pha}=\Cscr_X$ is known by Goldring-Koskivirta \cite{Goldring-Koskivirta-global-sections-compositio} \cite{Goldring-Koskivirta-rank-2-cones} gives:
\begin{corollary}
\label{cor-intro-pha-rank-2}
Assume that every $\fp$-simple factor $G_i$ of the adjoint group $G^{\ad}$ is either $\fp$-split of rank $2$ or that $G^{\ad}_k \cong \PGL(2)^m$ for some $m \geq 1$. 
Suppose that $\mu$ is component-wise maximal (see \ref{sec-component-wise-maximal}).
\begin{enumerate}
\item If the Cone Conjecture assumptions~\ref{conj-cone}\ref{item-component-wise-maximal}-\ref{item-pseudo-complete} hold, then every effective class $\eta \in \Tsf^1(Y)$ is strata-effective.  
\item In particular, in the Shimura variety setting~\ref{sec-intro-Hodge-type}, if $X=S_K$ is a Hilbert modular variety, a Siegel threefold, or a Picard modular surface at a split prime, then every effective class $\eta \in \Tsf^1(Y)$ is strata-effective.   
\end{enumerate}

\end{corollary}
\subsubsection{Examples of effective divisors that are not strata-effective} 
The Cone Conjecture is known to hold by Goldring-Koskivirta \cite{Goldring-Koskivirta-divisibility} in several cases where $\Cscr_{\pha} \subsetneqq \Cscr_{\GZip^{\mu}}$, including Picard modular surfaces at an inert prime and Siegel modular varieties of dimension $6$. 
This gives examples of effective divisors $\eta \in \Tsf^1(Y)$ that are not strata-effective.
Koskivirta-Imai \cite{Goldring-Imai-Koskivirta-weights} have classified the pairs $(G,\mu)$ such that $\Cscr_{\pha}= \Cscr_{\GZip^{\mu}}$. For example, $\Cscr_{\pha}= \Cscr_{\GZip^{\mu}}$ holds for all $n$ if $(\type(G), \type(L))=(\Bsf_n, \Bsf_{n-1})$ or $(\Csf_n, \Csf_{n-1})$. In these cases, combining the Koskivirta-Imai result, the Cone Conjecture and~\ref{th-pha-strata-eff} gives: 
\begin{conjecture}
\label{conj-divisors-Bn-Cn}
If $(\type(G), \type(L))=(\Bsf_n, \Bsf_{n-1})$ or $(\Csf_n, \Csf_{n-1})$, then every effective divisor $\eta \in \Tsf^1(Y)$ is strata-effective.    
\end{conjecture}
\subsection{Strata-effectivity of effective tautological classes}
\label{sec-intro-curves}

\begin{theorem}
\label{th-linear-eff}
    Suppose that the Zip stratification is linear. If $\type(G)=\Gsf_2$, assume that $p \geq 5$. For every Zip period map~\eqref{eq-intro-zip-period}, 
    every effective class $\eta \in \Tsf^*(X)$ is strata-effective.
\end{theorem}
Aided by the Planches of Bourbaki \cite{bourbaki-lie-4-6}, it is an easy computation~\ref{prop-linear-classification} that the Zip stratification is linear if and only if $(\type(G),\type(L))$ is one of the pairs in~\ref{th-grif-eff} other than $(\Dsf_n, \Dsf_{n-1})$. In the latter case, the stratification is almost linear but there are two strata of middle-length $|\Phi \setminus \Phi_I|/4$.
\subsubsection{Strata-effectivity of curves}
The question of strata-effectivity for curves can be seen as dual to strata-effectivity of divisors. We say that a tautological curve  $C$ is $w$-ordinary if the intersection $C \cap Y_w$ with the flag stratum $Y_w$ is open, dense in $C$. For such $C$, the intersection of a divisor of a partial Hasse invariant on $\overline{Y}_w$ with $C$ is a non-negative multiple of the class of the minimal stratum $Y_1$. Using this, we show: 
\begin{theorem}[\ref{th-gen-ord-hilb},~\ref{th-gen-ord-Sp4},~\ref{th-gen-ord-unip21-inert}]
 \label{th-intro-curves} Assume that $X$ is projective over $k$ and that $\zeta$ is smooth and surjective.
 Let $C \subset Y$ be a tautological $w_0$-ordinary curve in the flag space. If $\type(G)=\Asf_1^d, \Csf_2$ or if $G$ is non-split of type $\Asf_2$, then $[C] \in \Tsf^1(Y)$ is strata-effective. 
\end{theorem}
In particular~\ref{th-intro-curves} applies to toroidal compactifications $S_K^{\Sigma}$ of Hilbert modular varieties, Siegel threefolds and Picard modular surfaces at an inert prime.

In the $\Asf_2$-unitary case, partial Hasse invariants do not suffice. The nefness of the pullback of the Hodge line bundle is used to complete the argument. Using the partial Hasse invariants diagram in the $\Csf_2$-case we note that the Hodge vector bundle is not nef on the Siegel modular threefold (\ref{subsec-non-nef-A2}).

\section*{Acknowledgements}
We are grateful to Michel Brion, Jean-Stefan Koskivirta, Ben Moonen, Daniel Qin, Stefan Reppen, Torsten Wedhorn and Paul Ziegler for helpful discussions.

We both thank the Knut \& Alice Wallenberg Foundation for its support under grant KAW 2018.0356. W.G. thanks the Knut \& Alice Wallenberg Foundation for its support under Wallenberg Academy Fellow grant KAW 2019.0256 and grant KAW 2022.0308. W.G. thanks  the Swedish Research Council for its support under grant \"AR-NT-2020-04924.
\section{Notation \texorpdfstring{\&}{and} Background}
\label{sec-background}
Throughout, let $p$ be a prime and let $k$ be an algebraic closure of $\fp$.
\subsection{Chow rings}
\label{sec-chow-rings}
For a scheme $X$, denote by $\CHQs(X)$ the Chow ring of $X$ as in \cite[8.3]{Fulton-Intersection-theory-book}, tensored with $\QQ$ over $\ZZ$ (\ie with rational coefficients). For an algebraic group scheme $G$ acting on a scheme $X$, denote by $\CH_{G,\QQ}^{*}(X)$ the equivariant Chow ring (with rational coefficients) as in Edidin-Graham \cite[2.2-2.5]{Edidin-Graham-Equiv-Intersection-Theory}. For $G = \{e\}$ trivial and $X$ a scheme this recovers the usual Chow ring  $\CH_{G,\QQ}^{*}(X) = \CHQs(X)$. For a quotient stack $\Ycal$ admitting a presentation $\Ycal = [G\backslash X]$ with $X$ smooth (\ie $\Ycal$ is a smooth quotient stack), define $\CHQs(\Ycal) := \CH_{G,\QQ}^{*}(X)$. 
This is independent of the choice of smooth presentation  \cite[\Prop~16]{Edidin-Graham-Equiv-Intersection-Theory} 
\subsection{Sections and zero schemes}
\label{sec-sections-zero-scheme}
Let $X$ be a scheme, $\Lscr/X$ a line bundle and $s \in H^{0}(X,\Lscr)$ a global section. The section $s$ defines a morphism of sheaves $s : \mathcal{O}_{X}\rightarrow \mathcal{L}$. Taking duals gives a morphism of sheaves $s^{\vee} : \Lscr^{\vee} \rightarrow \Ocal_{X}$ and we define as usual $\Zero(s)\subset X$ to be the closed subscheme of $X$ defined by the ideal sheaf $\im(s^{\vee}) \subset \Ocal_{X}$. We say $s \in H^{0}(X,\Lcal)$ is injective if the induced map $s : \mathcal{O}_{X}\rightarrow \mathcal{L}$ is injective.  An injective section $s \in H^{0}(X,\mathcal{L})$ gives the relation $[\Zero(s)] = c_{1}(\mathcal{L}).[X]$ in the Chow ring $\CHQs(X)$. We often use this in the situation when $\iota : X \hookrightarrow X^{'}$ is a closed subscheme and $\Lscr = \iota^{*}\Lscr^{'}$ is the pullback of a line bundle on $X^{'}$. In this case we have, by pushfoward and projection formula, $\iota_{*}[\Zero(s)] = c_{1}(\Lscr^{'}).[X] \in \CHQs(X^{'})$.
 \subsection{Reductive groups, root data and Weyl groups}
 Let $G$ be a connected, reductive $\fp$-group. Let $\varphi: G \to G$ denote the relative Frobenius.
 
\subsubsection{Root data}
\label{sec-root-data}
Let $T$ be a maximal torus of $G$. Let $(\Xsf^*(T), \Phi; \Xsf_*(T), \Phi^{\vee})$ be the root datum of $(G_k,T_k)$, where $\Xsf^*(T)$ (resp. $\Xsf_*(T)$) is the character (resp. cocharacter) group of $T_k$ and $\Phi$ (resp. $\Phi^{\vee}$) is the  set of roots (resp. coroots) of $T_k$ in $G_k$.
Write $\langle,\rangle: \Xsf^*(T) \times \Xsf_*(T) \to \ZZ$ for the natural perfect pairing.
\subsubsection{Based root data}
\label{sec-based-root-data}
Let $B$ be a Borel subgroup of $G$ containing $T$. Let $\Phi^{+}\subseteq \Phi$ be the system of positive roots such that $\alpha \in \Phi^+$  if and only if the root group $U_{-\alpha}\subseteq B$. Let $\Delta \subset \Phi^+$ denote the base of simple roots. Let $\Xsf^*_+(T)$ be the cone of $\Delta$-dominant characters. 
\subsubsection{Derived subgroup and its simply-connected cover}
 \label{sec-derived-simply-connected}
  Let $G^{\der}$ be the derived subgroup of $G$, and $\tilde{G}$ the simply-connected cover of $G^{\der}$.  Let $\iota:\tilde G \to G$ be the composition of the isogeny $\tilde G \to G^{\der}$ with the inclusion $G^{\der} \hookrightarrow G$. For every subgroup $H \subset G$, set $\tilde H=\iota^{-1}H$ to be its preimage in $\tilde G$. If $H$ is a maximal torus (resp. Borel, parabolic or Levi subgroup) of $G$, then $\tilde H$ is a maximal torus (resp. Borel, parabolic, Levi subgroup) of $\tilde G$. By definition, $\tilde G$ is a group whose root datum satisfies  $\Xsf_*(\tilde T)=\ZZ\Phi^{\vee}(\tilde G, \tilde T)$.

\subsubsection{Weyl groups}
For $\alpha \in \Phi$, let $s_{\alpha} \in \aut(\Xsf^*(T)_{\QQ})$ be the root reflection $s_{\alpha}(x)=x-\langle x, \alpha^{\vee} \rangle \alpha$.   Let $W=W(G,T)$ be the Weyl group of $G_k$ relative $T_k$. Write $l:W \to \NN$ for the length function of the Coxeter system $(W, \{s_{\alpha}\}_{\alpha \in \Delta})$ and $w_0$ for its longest element.  
\subsubsection{Sub-root systems and parabolic subgroups}
Let $I \subset \Delta$. Set $\Phi_I:=\ZZ I \cap \Phi$. Let $W_I=\langle s_{\alpha} | \alpha \in I\rangle$ be the associated standard parabolic subgroup of $W$ and $w_{0,I}$ its longest element. Let ${}^{I}W$ denote set of minimal length representatives of the cosets  $W_I\backslash W$.  The longest element in ${}^{I}W$ is $w_{0,I}w_{0}$ of length $l(w_{0,I}w_{0})=(\Phi \setminus \Phi_I)/2$. The type of a standard parabolic subgroup $B \subset P \subset G$ is $\type(P)=\{\alpha \in \Delta | U_{\alpha} \subset P\}$.
\subsection{Special characters (and cocharacters)}
\label{sec-cond-characters}
The following definitions~\ref{def-orb-p-close} and~\ref{def-quasi-constant} were introduced in \cite[N.5.3]{Goldring-Koskivirta-Strata-Hasse} (see also \cite[1.1.2, (5.1.1)]{Goldring-Koskivirta-quasi-constant}).
\begin{definition}
\label{def-orb-p-close}
A cocharacter $\mu \in \Xsf_*(G)$ is \underline{orbitally $p$-close} if for some (equivalently every) maximal torus $T \subset G$ containing the image of $\mu$ and all roots $\alpha \in \Phi$ satisfying $\langle \alpha, \mu \rangle \neq 0$ one has 
\addtocounter{equation}{-1}
\begin{subequations}
\begin{equation}
\label{eq-def-orb-p-close}
    \frac{| \langle \sigma \alpha, \mu \rangle |}{| \langle \alpha, \mu \rangle |} \leq p-1 \textnormal{ for all } \sigma \in W \rtimes \gal(k/\fp)
\end{equation}        
\end{subequations}

\end{definition}
A character $\chi$ of a maximal torus $T$ of $G$ is orbitally $p$-close if~\eqref{eq-def-orb-p-close} holds for $\chi$ with coroots instead of roots.
\begin{definition}
\label{def-quasi-constant}
A cocharacter (resp. character) is \underline{quasi-constant} if it is orbitally $2$-close.
\end{definition}
 By \cite[1.4.4]{Goldring-Koskivirta-quasi-constant} the Hodge character of a Hodge-type Shimura variety is quasi-constant and the statement is independent of the choice of symplectic embedding. Note that the Hodge character is cominuscule but may fail to be minuscule -- this happens in the Siegel case.
 \subsubsection{Component-wise maximal cocharacters}
 \label{sec-component-wise-maximal}
 Following \cite{Goldring-Koskivirta-rank-2-cones}, a cocharacter $\mu \in \Xsf_*(G)$ is \underline{component-wise maximal} if, for every simple adjoint factor $G_i$ of $G_k^{\ad}$, either the projection $\mu_i$ of $\mu$ onto $G_i$ is central in $G_i$ or its centralizer $\cent_{G_i}(\mu_i)$ is a maximal Levi subgroup. Every minuscule, cominuscule and quasi-constant cocharacter is component-wise maximal. In particular, the conjugacy class of cocharacters associated to a Shimura datum by Deligne is component-wise maximal.
 \subsection{The stacks \texorpdfstring{$\GF^{\mu, P_0}$}{} of partial zip flags}
 \subsubsection{The stack $\GZip^{\mu}$, d'apr\`es Pink-Wedhorn-Ziegler \cite{Pink-Wedhorn-Ziegler-zip-data} \cite{Pink-Wedhorn-Ziegler-F-Zips-additional-structure}}
  \label{sec-partial-zip-flag}

 Let $P, P^+$ be the pair of opposite parabolic subgroups of $G$ of non-positive (resp. non-negative) $\mu$-weights. Let $L=P\cap P^+=\cent_{G_k} \mu$ be their common Levi factor.  As in \cite[1.2.3]{Goldring-Koskivirta-Strata-Hasse} assume that there exists an $\fp$-Borel subgroup $B \subset P$. Set $I=\type(P)$.
 Let $Q=(P^+)^{(p)}$. Recall  the zip group
 $$ E=\{(a,b) \in P \times Q | \varphi(\overline{a})=\overline{b}\},
 $$
 where $\overline{a}, \overline{b}$ denote the projections to the Levi factors. Then $\GZip^{\mu}:=[E \backslash G]$.
\subsubsection{The stack $\GF^{\mu.P_0}$, d'apr\'es Goldring-Koskivirta \cite{Goldring-Koskivirta-zip-flags}}
 Let $B \subset P_0 \subset P$ be an intermediate parabolic subgroup and set $I_0:=\type(P_0)$. The Zip group of $P_0$ is the subgroup $E_{P_0}=E \cap (P_0 \times Q)$ of $E$. Then 
 $$\GF^{\mu, P_0}:=[E_{P_0} \backslash G].
 $$
 When $P_0=B$, $I_0=\emptyset$ and $\GF^{\mu,P_0}=\GF^{\mu}$ is the stack of zip flags initially defined in \cite{Goldring-Koskivirta-Strata-Hasse}.
 The projections
 $$\GF^{\mu} \to \GF^{\mu,P_0} \to \GZip^{\mu}
 $$
 are fibrations with flag variety fibers $P_0/B$ and $P/P_0$ respectively. 
\subsubsection{Stratification}
\label{sec-fine-strat}
There is a Zip datum $\Zcal_0$ (usually not of cocharacter type) and a smooth surjective morphism
\addtocounter{equation}{-1}
\begin{subequations}
\begin{equation}
\label{eq-non-cocharacter}
 \GF^{\mu, P_0} \to \GZip^{\Zcal_0}   
\end{equation}    
\end{subequations} induced from the inclusion of zip groups $E_{P_0} \subset E_{\Zcal_0} = \{(x,y)\in P_0\times Q_0 \mid \varphi(\overline{x}) = \overline{y}\}$.
The zip stratification of the latter induces a stratification of the former parameterized by $\leftexp{I_0}{W}$, called the \underline{fine stratification} in \cite{Goldring-Koskivirta-zip-flags}. When $I_0=\emptyset$, the partial order is the Bruhat-Chevalley order, but when $I_0 \neq \emptyset$ it may be finer, as in the base case $I=I_0$.

\subsubsection{The compact dual}
\label{sec-cpt-dual}
Given a Frobenius zip datum $\Zcal=(G,P,L,Q,M, \varphi)$ (not necessarily of cocharacter type), its compact dual is the flag variety $\XX^{\vee}(\Zcal):=\GG_{\CC}/\PP_{\CC}$ over $\CC$, where $\GG_{\CC}$ is a connected reductive $\CC$-group of the same type as $G$ (\eg take $\GG$ adjoint of type $\type(G)$) and $\PP_{\CC}$ is a parabolic subgroup whose type is the same as that of $P$. If $\Zcal$ arises from $(G,\mu)$, write $\XX^{\vee}=\XX^{\vee}(G,\mu)$. 
\subsubsection{Associated bundles}
\label{sec-associated-bundles}
Write $L_{0}$ for the Levi of $P_{0}$. 
Given a $k$-representation $r$ of $P_{0}$, let $\Wscr(r)$ be the  ‘automorphic’ vector bundle on $\GF^{\mu,P_0}$ 
by applying the `associated sheaves construction' to the pullback of $r$ to $E_{P_0}$ via the projection $E_{P_{0}} \rightarrow P_{0}$. We consider three types of $P$-representations:
\begin{enumerate}
\item If  $r\in \Xsf^*(P_0)$ then write $\Lscr(r) := \Wscr(r)$ for the associated line bundle.
\item
\label{item-G-rep}
If $r$ is a $G$-representation, also write $\Wscr(r)=\Wscr(r_{|P_0})$ for the bundle associated to the restriction of $r$ to $P_0$.   
\item If $\lambda \in \Xsf^*(B)$ is $I_0$-dominant, let $\Vscr(\lambda)=\Wscr(\Ind_{B}^{P_0} \lambda)$.
\end{enumerate}
  
 Note that $\Ind_{B}^{P_0} \lambda$ is trivial on the unipotent radical $R_uP_0$ and its restriction to $L_0$ is $\Ind_{L_0 \cap B}^{L_0} \lambda$. In the Shimura variety setting~\ref{eq-zeta-shimura}, the pullbacks $\zeta^*\Wscr(r)$ of the bundles~\ref{item-G-rep} are the flat automorphic bundles. 
\subsubsection{Chow ring}
\label{sec-chow-zip-flag}
Let $[L_{\varphi} \backslash G]$ be the quotient stack of $G$ by $L$ acting by $\varphi$-conjugation.
Brokemper \cite{Brokemper-tautological} shows that 
\addtocounter{equation}{-1}
\begin{subequations}
\begin{equation}
\label{eq-brokember-iso}
\CHQs(\GZip^{\mu}) \cong \CHQs[L_{\varphi} \backslash G].    
\end{equation}    
He also shows that the classes of zip strata closures of codimension $i$ in $\GZip^{\mu}$ form a basis of $\CHQ^i(\GZip^{\mu})$.  Wedhorn-Ziegler \cite{Wedhorn-Ziegler-tautological} show that 
\begin{equation}
\label{eq-iso-L-conj-cpt-dual}    
\CHQs([L_{\varphi} \backslash G]) \cong H^{2*}(\XX^{\vee}(\CC), \CC).
\end{equation}
 Their arguments apply more generally to $\GZip^{\Zcal_0}$. Brokemper's argument is based on:
\end{subequations}

\begin{lemma}[{\cite[1.4.7]{Brokemper-tautological}}]
    \label{le-brok-ignore-unip}
    Let $G$ be a split extension of algebraic groups \[\begin{tikzcd}
        1 \arrow[r]& K \arrow[r]& G \arrow[r]& H \arrow[r]& 1
    \end{tikzcd}\] with $K$ unipotent. Let $X$ be a normal quasi-projective $G$-scheme. The map induced by a splitting $H \hookrightarrow G$ is an isomorphism $$\CH_{H,\QQ}^{*}(X) \stackrel{\sim}{\longrightarrow} \CH_{G,\QQ}^{*}(X).$$ 
\end{lemma}
\begin{corollary}
\label{cor-chow-zip-flag}
Via~\ref{le-brok-ignore-unip}, the projections $E_{P_0} \to L_0$ and $E_{\Zcal_0} \to L_0$ induce isomorphisms on Chow rings
$$\CHQs(\GZip^{\mu,P_0}) \cong \CHQs[L_{0,\varphi} \backslash G] \cong \CHQs(\GZip^{\Zcal_0}).$$    
\end{corollary}
One can also show that  $\CHQs(\GZip^{\Zcal_0}) \cong  H^{2*}(\XX^{\vee}(\Zcal_0)(\CC))$, where $\XX^{\vee}(\Zcal_0)=\GG_{\CC}/\PP_{0,\CC}$.   
\subsubsection{Stratification of smooth Zip period maps} 
\label{sec-strat-zip-period}
Write 
\addtocounter{equation}{-1}
\begin{subequations}
\begin{equation}
\Xcal:=\GZip^{\mu}, \quad \Xcal^{P_0}:=\GF^{\mu,P_0} \quad \textnormal{and} \quad \Ycal:=\Xcal^B:=\GF^{\mu}.    
\end{equation} For $w \in \iw$ (resp. $w \in \leftexp{I_0}{W}$, $w\in W$) write $\Xcal_w$ (resp. $\Xcal^{P_0}_w$, $\Ycal_w$) for the stratum parameterized by $w$.

For a Zip period map $\zeta:X \to \Xcal$~\eqref{eq-intro-zip-period}, let $X^{P_0}=X \times_{\Xcal} \Xcal^{P_0}$ and $Y=X \times_{\Xcal} \Ycal$ be the base changes of $X$ to the stacks of partial and full flags $\Xcal^{P_0}, \Ycal$ respectively. The spaces $X^{P_0}$ and $Y$ are called the partial and (full) flag spaces of $X$ respectively. Set $\zeta_{P_0}$ (resp. $\zeta_Y$) to be the base change of $\zeta$ to $X^{P_0}$ (resp. $Y$):
\begin{equation}
\xymatrix{
Y \ar[d] \ar[r]^{\zeta_Y} & \Ycal \ar[d] \\
X^{P_0} \ar[d] \ar[r]^{\zeta_{P_0}} & \Xcal^{P_0} \ar[d] \\
X  \ar[r]^{\zeta} & \Xcal
}    
\end{equation}    
\end{subequations}
Write $X_w:=\zeta^{-1}(\Xcal_w)$ (resp. $X^{P_0}_w:=\zeta_{P_0}^{-1}(\Xcal^{P_0}_w), Y_w=\zeta^{-1}(\Ycal_w)$) for the preimages of the strata. If $\zeta$ is smooth, then the $X_w$, (resp. $X^{P_0}_w$, $Y_w$) form stratifications of $X$ (resp. $X^{P_0}$, $Y$).
\subsection{Generalized Griffiths modules and bundles}
\label{sec-griffiths-bundle}
This section follows \cite[3.1]{Goldring-griffiths-bundle-group-theoretic}. 
\subsubsection{The Griffiths module of a filtered vector space}
Let $F$ be a field, let $V$ be a finite-dimensional $F$-vector space and let $\fil^{\bullet}$ be a finitely supported, descending filtration on $V$ indexed by $\ZZ$. The Griffiths module of $(V, \fil^{\bullet})$ is 
\begin{equation}
\Grif(V, \fil^{\bullet}):=\sum_{n \in \ZZ} \fil^i.    
\end{equation}
\subsubsection{Deligne's convention}
\label{sec-deligne-convention}
Following Deligne's sign convention for Hodge structures \cite{Deligne-Shimura-varieties}, a cocharacter 
$\mu \in \Xsf_*(\GL(V))$ is assigned the filtration $\fil^{\bullet}$ where $\fil^{i}$ is the $-i$-weight space of $\mu$.
\subsubsection{The Griffiths module of a triple \texorpdfstring{$(G,\mu,r)$}{!}}
\label{sec-grif-module-triple}
Let $G$ be a connected, reductive $F$-group, $\mu \in \Xsf_*(G)$ and $r:G \to \GL(V)$ an $F$-representation of $G$. 
The Griffiths module $\Grif(G,\mu,r)$ is the Griffiths module of $(V, \fil^{\bullet})$, where $\fil^{\bullet}$ is the filtration associated to $r \circ \mu \in \Xsf_*(\GL(V))$ by the rule~\ref{sec-deligne-convention}. 
Let $L:=\cent(\mu)$. 
Then $\Grif(G, \mu, r)$ is an $L$-module. 
\subsubsection{The Griffiths character of $(G,\mu,r)$} The Griffiths character is the determinant of the Griffiths module: $\grif(G,\mu,r)=\det\Grif(G,\mu,r) \in \Xsf^*(L)$.
\subsubsection{Griffiths bundles on $\GZip^{\mu}$}
\label{sec-grif-GZip}
Apply~\ref{sec-grif-module-triple} with $F=k$ algebraically closed.
The Griffiths-Hodge vector bundle $\Grif(\GZip^{\mu},r)$ (resp. line bundle $\grif(\GZip^{\mu},r)$) of $\GZip^{\mu}$ relative $r$ is the vector bundle (resp. line bundle) on $\GZip^{\mu}=[E\backslash G]$ associated as in~\ref{sec-associated-bundles} to the $L$-module $\Grif(G,\mu,r)$ (resp. character $\grif(G,\mu,r)$ of $L$).
\subsubsection{Griffiths bundles of Zip period maps}
\label{sec-def-grif-zip-period}
By definition, the Griffiths bundle  of a Zip period map~\ref{eq-intro-zip-period} relative $r$ is the pullback $$\Grif(X,\zeta,r):=\zeta^*\Grif(\GZip^{\mu},r).$$
\subsection{The Schubert stack and Chevalley's formula}
\label{sec-chevalley-formula}
Let $\Sbt:=[B \backslash G/B]$. Denote $\Lscr_{\Sbt}(\lambda_1,\lambda_2)$ the line bundle on $\Sbt$ induced from a pair of characters $(\lambda_1,\lambda_2) \in \Xsf^*(T)\times \Xsf^*(T)$.
Denote the set of lower neighbours $E_{w} := \{\alpha \in \Phi^{+} \mid ws_{\alpha} < w \text{ and } l(ws_{\alpha}) = l(w) - 1\}$. For the sake of exposition denote \begin{equation*}
    \Mscr_{w}(\lambda) := \Lscr_{\Sbt}(\lambda,-w^{-1}\lambda)
\end{equation*} 
\begin{lemma}[Chevalley's formula {\cite[\Prop~10]{Chevalley-decompositions-cellulaires}, \cf \cite[2.2.1]{Goldring-Koskivirta-Strata-Hasse} for the stack $\Sbt$}]
    \label{le-chevalley-form}
    Consider $s \in H^{0}(\Sbt_{w},\Mscr_{w}(\lambda))$ as a rational section on $\overline{\Sbt}_{w}$.
    Then $$div(s) = -\sum_{\alpha \in E_{w}}\langle\lambda,w\alpha^{\vee}\rangle[\overline{\Sbt}_{ws_{\alpha}}].$$ 
\end{lemma}
\begin{lemma}
\label{le-admitpha}
    Let $w \in W$ and $\lambda \in \CharT$. Then $H^{0}(\overline{\Sbt}_{w},\Mscr_{w}(\lambda)) \neq 0$ if and only if $\langle\lambda,w\alpha^{\vee}\rangle \leq 0$ for all $\alpha \in E_{w}$.
\end{lemma}
\subsubsection{Proof}. A global section on $\Sbt_{w}$ gives a rational section $s$ on $\overline{\Sbt}_{w}$ whose divisor is given by the Chevalley formula \ref{le-chevalley-form}. If there is a nonzero global section $s \in H^{0}(\Sbt_{w},\Mscr_{w}(\lambda))$ on $\overline{\Sbt}_{w}$ then $\div(s)$ is effective and all coefficients in the Chevalley formula are non-negative. Conversely, if the coefficients are all non-negative then the Cartier divisor given by the rational section $s$ is effective and thus cut out by a global section on $\overline{\Sbt}_{w}$. So $H^{0}(\overline{\Sbt}_{w},\Mscr_{w}(\lambda)) \neq 0$.     
\qed

We refer to $f_{\lambda,w} \in H^{0}(\overline{\Sbt}_{w},\Mscr_{w}(\lambda))$ as Chevalley or highest weight sections. The Chevalley formula gives the following expression for the product $c_{1}(\Mscr_{w}(\lambda)).[\overline{\Sbt}_{w}]$ in terms of the classes of the divisors $\overline{\Sbt}_{ws_{\alpha}}$ of $\overline{\Sbt}_{w}$.
\begin{lemma}
\label{le-cyc-cl-sbt}
    Let $w\in W$ and $\lambda \in \CharT$. $$c_{1}(\Mscr_{w}(\lambda)).[\overline{\Sbt}_{w}] = -\sum_{\alpha \in E_{w}}\langle\lambda,w\alpha^{\vee}\rangle[\overline{\Sbt}_{ws_{\alpha}}]$$ in $\CHQ^*(\Sbt)$.
\end{lemma}
\subsubsection{Proof}
    Let $\chi \in X^*_+(T)$ be dominant and regular, \eg $\chi=\sum_{\alpha \in \Phi^+}\alpha$. Then for some $m>0$, the character $\lambda + m\chi$ is also dominant and regular. Line bundles given by dominant and regular characters admit Chevalley-highest weight sections which are injective and so there exist injective sections $s \in H^{0}(\overline{\Sbt}_{w},\Mscr_{w}(m\chi))$ and $t \in H^{0}(\overline{\Sbt}_{w},\Mscr_{w}(\lambda + m\chi))$. Then 
    \begin{equation*}
        \begin{aligned}
        c_{1}(\Mscr_{w}(\lambda)).[\overline{\Sbt}_{w}] &= c_{1}(\Mscr_{w}(\lambda + m\chi)).[\overline{\Sbt}_{w}] - c_{1}(\Mscr_{w}(m\chi)).[\overline{\Sbt}_{w}] \\&= [\Zero(t)] - [\Zero(s)] \\&= -\sum_{\alpha \in E_{w}}\langle\lambda + m\chi,w\alpha^{\vee}\rangle[\overline{\Sbt}_{ws_{\alpha}}] + \sum_{\alpha \in E_{w}}\langle m\chi,w\alpha^{\vee}\rangle[\overline{\Sbt}_{ws_{\alpha}}] \\&= -\sum_{\alpha \in E_{w}}\langle\lambda,w\alpha^{\vee}\rangle[\overline{\Sbt}_{ws_{\alpha}}]
    \end{aligned}
    \end{equation*}
\qed
\subsection{Partial Hasse invariants and their cones}
\label{sec-pha}
Following \cite{Goldring-Koskivirta-rank-2-cones}, \cite{Goldring-Koskivirta-divisibility} \cite{Goldring-Koskivirta-GS-cone}, the "partial Hasse invariant" map is the smooth morphism
\begin{equation}
    \hsf:\GF^{\mu} \rightarrow \Sbt
\end{equation}
 constructed in \cite[(2.3.1)]{Goldring-Koskivirta-Strata-Hasse} and denoted $\psi$ there. 
\subsubsection{Partial Hasse invariants}
\label{sec-partial-Hasse}
\addtocounter{equation}{-1}
\begin{subequations}
A \underline{partial Hasse invariant} for $w \in W$ is a global section $\hsf^*f_{\lambda,w}$ over the flag stratum closure $\overline{\Ycal}_w$ which is the pullback of a Chevalley-highest weight section $f_{\lambda, w} \in H^{0}(\overline{\Sbt}_{w},\Mscr_{w}(\lambda))$ over a Schubert stratum closure $\overline{\Sbt}_w  \subset \Sbt$.

Define the partial Hasse invariant cone of $w \in W$ to be 
\begin{equation}
 \label{eq-def-pha-cone}   
 \Cscr_{\pha, w}:=\{\lambda \in \Xsf^*(T) \mid \Lscr(n\lambda)|_{\overline{\Ycal}_{w}} \text{ admits a partial Hasse invariant for some }n>0\}
\end{equation}  By the explicit description of $\hsf^*: \Pic(\Sbt) \rightarrow \Pic(\GF^{\mu})$ (\cite[\Lem~3.1.1(b)]{Goldring-Koskivirta-Strata-Hasse}), 
\begin{equation}
\label{eq-explicit-formula-Z-Chev}
\hsf^{*}\Mscr_{w}(\lambda) = \Lscr(D_{w}(\lambda))\hspace{0.5cm} \text{ and } \hspace{0.5cm}\Cscr_{\pha,w}=D_w(\Cscr_{\Sbt, w})
\end{equation} 
where $ D_{w}:\Xsf^*(T) \rightarrow \Xsf^*(T)$ is given by 
$\lambda \mapsto \lambda-p\leftexp{\sigma^{-1}}{(zw^{-1}\lambda)}$ and 
\begin{equation}
    \Cscr_{\Sbt,w} := \{\lambda \in \Xsf^*(T) \mid H^{0}(\overline{\Sbt}_{w},\Mscr_{w}(n\lambda))\neq 0  \text{ for some }n>0\}
\end{equation}
Given $\lambda \in \CharT$ there is some $n>0$ and $\chi \in \CharT$ such that $D_{w}(\chi) = n\lambda$ by \cite[3.1.3(a)]{Goldring-Koskivirta-Strata-Hasse}.

Lemma \ref{le-admitpha} gives the following condition for $\lambda$ to lie in the cone $\Cscr_{\pha,w}$: $\lambda \in \Cscr_{\pha,w}$ if and only $\lambda=\hsf^{-1}(\chi)$ and there is some $m>0$ such that $\langle m\chi,w\alpha^{\vee}\rangle \leq 0$ for all $\alpha \in E_{w}$.
\end{subequations}
\subsubsection{Partial Hasse invariants for partial flag spaces}
\label{sec-partial-hasse-partial-flag}
In the setting of~\ref{sec-partial-zip-flag}, let $w \in \leftexp{I_0}{W} $ and $\lambda \in \Xsf^*(L)$. A partial Hasse invariant for $(\lambda,w)$ over the flag stratum closure $\overline{\Xcal}^{P_0}_w$ is a global section $h \in H^0(\overline{\Xcal}^{P_0}_w, \Lscr(\lambda))$ whose pullback to $H^0(\overline{\Ycal}_w, \Lscr(\lambda))$ is a partial Hasse invariant. Let
\addtocounter{equation}{-1}
\begin{subequations}
 \begin{equation}
 \label{eq-def-partial-pha-cone}   
 \Cscr_{\pha, w}^{P_0}:=\{\lambda \in \Xsf^*(T) \mid \Lscr(n\lambda)|_{\overline{\Xcal}^{P_0}_{w}} \text{ admits a partial Hasse invariant for some }n>0\}
\end{equation}
\end{subequations}
\subsubsection{
\texorpdfstring{$\textdbend$}{!}}  The cones $\Cscr_{\pha,w}$ were previously called Schubert cones and denoted $\Cscr_{\Sbt, w}$ in \cite{Goldring-Koskivirta-global-sections-compositio,Koskivirta-automforms-GZip}.
Note that there is a sign error in \cite{Goldring-Koskivirta-global-sections-compositio}, where $\sigma$ appears instead of $\sigma^{-1}$; this error does not affect \cite{Goldring-Koskivirta-global-sections-compositio} because $\sigma$ has order dividing two in the cases where~\eqref{eq-explicit-formula-Z-Chev} is applied there.

\subsection{The Cone Conjecture, d'apr\texorpdfstring{\`}{}es Goldring-Koskivirta \texorpdfstring{\cite{Goldring-Koskivirta-global-sections-compositio}}{}}
Recall that a $k$-stack $S$ is pseudo-complete if every global function $f \in H^0(S,\Ocal_S)$ is locally constant. 
\begin{conjecture}[{\cite[2.1.6]{Goldring-Koskivirta-global-sections-compositio},~\cite[3.2.1]{Goldring-Koskivirta-rank-2-cones}}]
\label{conj-cone}
    Let $\zeta : X \rightarrow \GZip^\mu$ be a Zip period map~\eqref{eq-intro-zip-period}. Suppose that 
\begin{enumerate}
\item
\label{item-component-wise-maximal}
The cocharacter $\mu$ is component-wise maximal~\eqref{sec-component-wise-maximal};
\item 
\label{item-smooth}    
The morphism $\zeta$ is smooth and its restriction to every connected component $X^+\subset X$ is surjective;
\item
\label{item-pseudo-complete}
    For every simple root $\alpha \in \Delta$, the stratum closure $\overline{Y}_{s_{\alpha}}$ is pseudo-complete.
\end{enumerate}
Then $\Cscr_{X} = \Cscr_{\GZip^\mu}$.
\end{conjecture}
\begin{rmk} \
\begin{enumerate}
\item The pseudo-completeness condition~\ref{item-pseudo-complete} automatically holds if $X$ (hence also $Y$) is proper over $k$. 
In particular, by~\ref{sec-intro-Hodge-type} the hypotheses are satisfied if $X=S_K^{\Sigma}$ is a sufficiently nice toroidal compactification of a Hodge-type Shimura variety.
\item In \cite{Goldring-Imai-Koskivirta-weights} it is determined when $\Cscr_{\GZip^\mu} = \Cscr_{\pha}$ in terms of $G$ and $\mu$.

\end{enumerate}
    
\end{rmk}

\section{Functoriality}
\label{sec-functoriality}
\subsection{The Wedhorn-Ziegler isomorphism}
Observe that the proof of Wedhorn-Ziegler actually shows that the isomorphism~\ref{eq-Wedhorn-Ziegler-iso} holds for more general Zip period maps:
\begin{theorem} Let $B \subset P_0 \subset P$ be an intermediate parabolic.  Assume that: 
\begin{enumerate}
\item The morphism $\zeta: X \to \GZip^{\mu}$ satisfies~\ref{conj-cone}\ref{item-smooth}
\item
\label{item-hyp-minimal-nonzero}
The minimal stratum $X_e \subset X$ has nonzero class $[X_e] \neq 0$ in $\Tsf^{|\Phi \setminus \Phi_I|/2}(X)$.
\end{enumerate}
  Then pullback by $\zeta_{P_0}$~\eqref{sec-strat-zip-period} is a ring isomorphism onto the tautological ring:
\addtocounter{equation}{-1}
\begin{subequations}
\begin{equation}
\label{eq-wedhorn-ziegler-generalized}
    \zeta_{P_0}^*:\CHQs(\GF^{\mu,P_0}) \stackrel{\sim}{\longrightarrow} \Tsf^*(X^{P_0}).   
    \end{equation}    
\end{subequations}  
   \end{theorem}
\begin{rmk}
In particular hypothesis~\ref{item-hyp-minimal-nonzero} holds if $X$ is a proper $k$-scheme.    
\end{rmk}
Recall the full flag space  $Y=X^P$. The extremal cases $B=P_0$ and $P_0=P$ give the pullback isomorphisms:$$\zeta^*:\CHQs(\GZip^{\mu}) \stackrel{\sim}{\longrightarrow} \Tsf^*(X) \textnormal{ and } \zeta_Y^*\CHQs(\GZip^{\mu}) \stackrel{\sim}{\longrightarrow} \Tsf^*(Y).$$ 
\subsubsection{Proof of~\ref{eq-wedhorn-ziegler-generalized}}
The proof of Wedhorn-Ziegler that~\eqref{eq-wedhorn-ziegler-generalized} holds when $P_0=P$ and $X=S_K^{\Sigma}$ is the special fiber of a smooth toroidal compactification a Shimura variety is comprised of three parts:
\begin{enumerate}
\item
\label{item-wz-cpt-dual}
The isomorphism with the cohomology of the compact dual $\CHQs(\GZip^{\mu}) \cong H^{2*}(\GG_{\CC}/\PP_{\CC}, \QQ)$
\item
\label{item-wz-reduce-minimal-nonzero}
The observation that, due to~\ref{item-wz-cpt-dual}, it suffices to show that the pullback of the class $[\Xcal_e]$ of the minimal stratum $\Xcal_e \subset \Xcal=\GZip^{\mu}$ is nonzero.
\item
\label{item-wz-minimal-nonzero}
Since $\zeta$ is smooth $\zeta^*[\Xcal_e]=[S_{K,e}^{\Sigma}]$ and $[S_{K,e}^{\Sigma}] \neq 0$ since $S_K^{\Sigma}$ is proper over $k$.
\end{enumerate}
In fact, in~\ref{item-wz-cpt-dual}, the isomorphism~\eqref{eq-iso-L-conj-cpt-dual} suffices for~\ref{item-wz-reduce-minimal-nonzero}. 
We observed~\ref{cor-chow-zip-flag} that the isomorphism~\eqref{eq-iso-L-conj-cpt-dual} generalizes to the stacks $\GF^{\mu,P_0}$. 
The hypthesis that $[X_e] \neq 0$ implies that also $[X^P_e] \neq 0$, since the restriction of the projection $X^{P_0} \to X$  to $X^P_e$ is finite onto $X_e$. 
Thus~\ref{eq-wedhorn-ziegler-generalized} follows from Wedhorn-Ziegler's argument due to the hypotheses that $\zeta$ is smooth and that $[X_e] \neq 0$. \qed 
\begin{rmk} In \cite[Proof of 7.12]{Wedhorn-Ziegler-tautological} there is a more complicated argument to show that $\zeta^*[\Xcal_e]$ is nonzero. It seems that this argument was meant to prove that $\zeta^*$ was injective under hypotheses that were weaker than the smoothness of $\zeta$, as in an earlier version of \cite{Wedhorn-Ziegler-tautological} which preceded Andreatta's proof that $\zeta^{\Sigma}$ was smooth for a smooth toroidal compactification $S_K^{\Sigma}$~\eqref{sec-intro-Hodge-type}.
    
\end{rmk}
\subsection{Products I: \texorpdfstring{$\GZip^{\mu}$}{G-Zip stacks} and tautological rings}
\label{sec-gen-products}
Assume that $G=G_1 \times G_2$ is a direct product of two connected, reductive $\fp$-groups.
\subsubsection{Product decomposition of \texorpdfstring{$\GZip^{\mu}$}{}}
Let $\mu_i$ be the projection of $\mu$ onto $G_i$. Let $\iw_1$ and $\iw_2$ be the corresponding sets of minimal representatives.  It is straightforward (\cf \cite[5.2.3]{Goldring-Koskivirta-rank-2-cones}) that $\iw=\iw_1 \times \iw_2$ and that the stack of $G$-Zips decomposes as a product of the corresponding stacks of $G$-Zips: 
\addtocounter{equation}{-1}
\begin{subequations}
\begin{equation}
\label{eq-product-GZipmu}
\GZip^{\mu}=\GoneZip^{\mu_1}\times \GtwoZip^{\mu_2}.  
\end{equation}    
\end{subequations}

\begin{lemma}
\label{lem-Chow-GZip-prod}
The exterior product induces a Kunneth isomorphism of graded rings
\addtocounter{equation}{-1}
\begin{subequations}
\begin{equation}
\label{eq-product-chow-GZipmu} 
\CHQs(\GoneZip^{\mu_1}) \otimes \CHQs(\GtwoZip^{\mu_2})
\stackrel{\sim}{\longrightarrow} 
\CHQs(\GZip^{\mu}), \end{equation}    
\end{subequations}
which maps the basis $([\Xcal_v] \otimes [\Xcal_w])_{(v,w) \in \iw_1 \times \iw_2}$ to the basis $(\Xcal_{v}\times \Xcal_{w})_{(v,w) \in \iw_1 \times \iw_2}$.
\end{lemma}
\subsubsection{Proof}
By \cite[4.12]{Wedhorn-Ziegler-tautological}, the Chow ring of $\GZip^{\mu}$ is isomorphic as a graded ring to the Betti cohomology of the compact dual: $\CHQs(\GZip^{\mu}) \cong H^{2*}(\XX^{\vee}(G,\mu)(\CC), \QQ)$. The cellular decomposition of flag varieties implies that the analogue of~\eqref{eq-product-chow-GZipmu} holds when the Chow rings of the stacks of $G$-Zips are replaced by the cohomology rings of the compact duals  (\cf Fulton \cite[\Exa~1.10.2]{Fulton-Intersection-theory-book}). 
\qed

An immediate corollary of~\ref{lem-Chow-GZip-prod} and the Wedhorn-Ziegler isomorphism~\ref{eq-Wedhorn-Ziegler-iso} is:
\begin{corollary} Assume that $\zeta:  X \to \GZip^{\mu}$ is smooth and surjective. Define graded subalgebras $T_i:=\zeta^*\CHQs(\GiZip^{\mu_i})   \otimes 1$ of the tautological ring $\Tsf^*(X)$. Then the tautological ring of $X$ is also a tensor product $$\Tsf^*(X) = T_1 \otimes T_2.$$
\end{corollary}
\subsubsection{Remark}
\label{rmk-problem-eff-product}
Suppose that for all smooth surjective $\zeta_i: X_i \to \GiZip^{\mu_i}$, $i=1,2$,  all effective tautological classes $\alpha_i \in \Tsf^*(X_i)$ are strata-effective. In general, we are unable to conclude from this hypothesis that, for every smooth surjective $\zeta:X \to \GZip^{\mu}$, all effective tautological classes $\alpha \in \Tsf^*(X)$ are strata-effective. The reason is that we have no way to decompose $\alpha \in \Tsf^*(X)$ \textit{effectively} based on the knowledge that $\GZip^{\mu}$ is a product~\eqref{eq-product-GZipmu}. However, we do give a criterion for showing that every effective tautological class $\alpha \in \Tsf^*(X)$ is strata-effective and if the criterion applies to $(G_i,\mu_i)$ then it also applies to $(G,\mu)$. In this way, strata-effectivity is propagated to products in all cases where it is shown to hold for the constituents. 
\subsection{Products II: Hodge-Chern classes}
\label{sec-prod-Hodge-Chern}
In contrast to the problem~\ref{rmk-problem-eff-product}, the Griffiths bundle is compatible with direct sums. This is the basis of:
\begin{proposition}
 \label{prop-Grif-product}
 Let $r_i:G_i \to \GL(V_i)$ be a $k$-representation of $G_i$. Let $V:=V_1 \oplus V_2$ and let $r:G \to \GL(V)$ be the sum of $r_1$ and $r_2$. If the Hodge-Chern classes of $\Grif(\GiZip^{\mu_i},r_i)$ are strata-effective for $i=1,2$, then so are the Hodge-Chern classes of $\Grif(G,\mu,r)$.
\end{proposition}
\subsubsection{Proof}
By the definitions~\ref{sec-griffiths-bundle}  $\Grif(G,\mu,r)=\Grif(G,\mu_1,r_1) \oplus \Grif(G,\mu_1,r_2)$ as $L$-modules and the Griffiths bundle is a sum $$\Grif(\GZip^{\mu},r)=\Grif(\GZip^{\mu}, r_1) \oplus \Grif(\GZip^{\mu}, r_2).$$
Moreover, $\Grif(\GZip^{\mu}, r_i)=\pi_{i}^*\Grif(\GiZip^{\mu_i}, r_i)$ where the $\pi_i$ are the projections onto the factors. So if the Chern classes of $\Grif(\GiZip^{\mu_i}, r_i)$ are strata-effective then the Chern classes of $\Grif(\GZip^{\mu},r_i)$ are supported on the strata closures $X_w \otimes 1$ (resp. $1 \otimes X_w$) with $w \in \iw_1$ (resp. $w \in \iw_2$).
Hence the result follows from the standard formula for the Chern classes of a direct sum and~\ref{lem-Chow-GZip-prod}. \qed 
\subsection{Equi-adjoint groups}
\label{sec-central-quotient}
\subsubsection{Adjoint morphisms}
Let $\psi:G \to H$ be a morphism of connected, reductive $\fp$-groups. If $\psi$ has central kernel, then $\psi^{-1}(\Zsf_H) \subset \Zsf_G$ (\cf \cite[2.2.16]{Goldring-propagating-algebraicity-functoriality}) Composition with the adjoint projection $H \twoheadrightarrow H^{\ad}$ induces the adjoint morphism $\psi^{\ad}:G/\psi^{-1}(\Zsf_H) \to H^{\ad}$. One says that $\psi$ induces an isomorphism on adjoint groups if $\psi^{\ad}$ is an isomorphism. By \cite[4.23]{Wedhorn-Ziegler-tautological}:
\begin{proposition}
\label{prop-central-quotient}
Let $\mu \in \Xsf_*(G)$ and let $P_0$ be an intermediate parabolic $B \subset P_0 \subset P$.
Assume that $\psi: G \to H$ induces an isomorphism of adjoint groups. Then 
\addtocounter{equation}{-1}
\begin{subequations}
\begin{equation}
\label{eq-equi-adjoint-pullback}    
\psi^*: \CHQs(\HF^{\psi_*\mu, \psi (P_0)}) \stackrel{\sim}{\longrightarrow} \CHQs(\GF^{\mu, P_0})
\end{equation}    
\end{subequations}
is an isomorphism of graded $\QQ$-algebras which maps the class of the $w$-stratum closure in $\HF^{\psi_*\mu, \psi(P_0)}$ to a rational multiple $>0$ of the class of the $w$-stratum closure in $\GF^{\mu, P_0}$. 
\end{proposition}
\subsubsection{Proof}
The case $P_0=P$ is \cite[4.23]{Wedhorn-Ziegler-tautological}. The same argument shows the analogue of~\eqref{eq-equi-adjoint-pullback} for the stacks $\GZip^{\Zcal_0}$~\eqref{sec-fine-strat}. It remains to observe that $\GF^{\mu, P_0} \to \GZip^{\Zcal_0}$ induces an isomorphism on Chow rings. This follows from Brokemper's Lemma~\ref{le-brok-ignore-unip}.
\qed
\subsection{Functoriality of the Griffiths bundle}
\label{sec-funct-grif}
The following is immediate from the definitions~\ref{sec-griffiths-bundle}; it explains that the Griffiths modules and bundles essentially depend only on the image of the representation $r$. 
\begin{lemma}
\label{lem-funct-grif}
Let $\psi: G \to H$ be a morphism of connected, reductive $k$-groups, $\mu \in \Xsf_*(G)$ and let $r:H_k \to \GL(V)$ be an $H_k$-module. Then 
$$ \Grif(G, \mu,r \circ \psi)=\psi^*\Grif(H, \psi\circ \mu, r) \textnormal{ and } \Grif(\GZip^{\mu},r \circ \psi)=\psi^*\Grif(\HZip^{\psi \circ \mu}, r).
$$
\end{lemma}
\subsection{Chern classes associated to \texorpdfstring{$\tilde G$}{tilde G}-representations}
\label{sec-tilde-G-reps}
Let $\iota:\tilde G \to G$ be the natural map~\eqref{sec-derived-simply-connected}. There exists $\tilde \mu \in \Xsf_*(\tilde G)$ such that $\GZip^{\iota_*\tilde\mu} = \GZip^{\mu}$. 
Let $\tilde P=\iota^*P$ and let $r:\tilde P \to \GL(V)$ be a $k$-representation. 
By~\ref{prop-central-quotient},for every $i$ there is a unique class in  $\CHQ^i(\GZip^{\mu})$ whose pullback to $\CHQ^i(\tildeGZip^{\tilde \mu})$ is $\csf_i \Wscr(r)$. 
In this way, one has Chern classes on $\GZip^{\mu}$ associated to a representation $r$ of $\tilde P$, even if it does not factor through $P$. In particular one has Chern classes $\csf_i\Grif(G,\mu,r)$ of the Griffiths bundle for every representation $r$ of $\tilde G$, even if $r$ is not the pullback to $\tilde G$ of a $G$-representation.     
\section{
Hodge-Chern classes I: Groups of type \texorpdfstring{$\Asf_1^m$}{a power of A1} \texorpdfstring{\&}{and} Hilbert modular varieties}
\label{sec-hilbert-strata-eff}
\S\ref{sec-hodge-chern-A1} explains how the isomorphisms between the Chow rings of $\GZip^{\mu}$ for equi-adjoint groups~\eqref{prop-central-quotient} and the Chern classes associated to $\tilde G$-representations~(\S\ref{sec-tilde-G-reps}) uniformly gives Hodge-Chern classes on arbitrary groups of type $\Asf_1^m$. The strata-effectivity of these classes is stated in this generality in \S\ref{sec-hilb-strat-eff}. \S\ref{sec-hilbert-case} explains how the Hilbert modular varieties of the form $S_K$~\ref{sec-intro-Hodge-type} and their smooth toroidal compactifications $S_K^{\Sigma}$ are a special case. The main strata-effectivity result~\ref{th-precise-eff-hilbert} is translated into a commutative algebra statement in \S\ref{sec-comm-alg}. Using \S\ref{sec-comm-alg}, the proof of~\ref{th-precise-eff-hilbert} given in \S\S\ref{sec-reciprocity}-\ref{sec-PI-in-terms-of-LJ} can be read independently of the rest of the paper. Some geometric motivation for two of the algebraic constructions in \S\S\ref{sec-reciprocity}-\ref{sec-PI-in-terms-of-LJ} is given in~\S\ref{sec-geometric}
\subsection{Hodge-Chern monomials for groups of type \texorpdfstring{$\Asf_1$}{A1}}
\label{sec-hodge-chern-A1}
Let $G$ be an $\fp$-group of type $\Asf_1^d$ for some $d \geq 1$. Assume that $\mu \in \Xsf_*(G)$ is dominant and regular for a choice of $\fp$-Borel pair $(B,T)$ such that $T_k$ contains the image of $\mu$. Let $\tilde \mu$ as in~\S\ref{sec-tilde-G-reps}.
\subsubsection{A Based root datum of $\tilde G$}
Identify $\tilde{G}_k \stackrel{\sim}{\longrightarrow} \SL(2)_k^{\ZZ/d}$ (\ie index the $d$ factors by $\ZZ/d$). Let $B,T$, $\iota: \tilde G \to G$, $\tilde B:=\iota^*B$ and $\tilde T :=\iota^*T$ as in \ref{sec-root-data}-\ref{sec-derived-simply-connected}. Then $\tilde B$ is an $\fp$-Borel subgroup of $\tilde G$ and $\tilde T$ is an $\fp$-maximal torus of $\tilde B$. Identify the based root datum of $(\tilde G_k,\tilde B_k,\tilde T_k)$ with $$(\ZZ^{\ZZ/d}, \{2e_i | i \in \ZZ/d\}; \ZZ^{\ZZ/d}, \{e_i | i \in \ZZ/d\}).$$ 
Without loss of generality, we may assume that $\tilde \mu=(1, \ldots ,1)$. 
\subsubsection{Degree one Hodge-Chern classes}
\label{sec-deg-1-Hodge-Chern}
By \S\ref{sec-tilde-G-reps}, $-e_i \in \Xsf^*(\tilde B)$ gives a well-defined class $l_i=\csf_1\Lscr(-e_i) \in \CHQ^1(\GZip^{\mu})$. For every $i \in \ZZ/d$, let $r_i$ be the two-dimensional $k$-representation of $\tilde G$ given as the composition of the projection onto the $i$th $\SL(2)_k$-factor with the inclusion into $\GL(2)_k$. Then the Griffiths line bundle $\grif(G, \mu,r_i)$ has first Chern class $l_i$ and all the higher Chern classes of $\Grif(G,\mu,r_i)$ vanish. 
\subsubsection{The Wedhorn-Ziegler isomorphism for type $\Asf_1^d$}
\addtocounter{equation}{-1}
\begin{subequations}
Let
\begin{equation}
\label{eq-taut-hilbert}
R := \QQ[l_0, \ldots, l_{d-1}]/(l_0^2, \ldots ,l_{d-1}^2)
\end{equation}
An easy special case of the Wedhorn-Ziegler isomorphism~\eqref{eq-Wedhorn-Ziegler-iso} gives an identification of graded $\QQ$-algebras
\begin{equation}
\label{eq-wedhorn-ziegler-hilbert}
  R \cong \CHQs(\GZip^{\mu}) \cong H^{2*}((\poc)^d, \QQ)
\end{equation}
\end{subequations}
Recall that~\eqref{eq-wedhorn-ziegler-hilbert} also follows more elementarily way using partial Hasse invariants, see \cite[\Rmk~3.2]{Cooper-tautological-ring-hilbert} and \ref{sec-PI-geometric}. 

\subsubsection{Hodge-Chern monomials}
\label{sec-hodge-chern-monomials}
\addtocounter{equation}{-1}
\begin{subequations}
Let $I \subset \ZZ/(d)$. Define $L_I\in \Gr^{|I|}R=\CHQ^{|I|}(\GZip^{\mu})$ by $L_{\emptyset}=1$ and   
\begin{equation}
L_I:=\prod_{i \in I}l_i \textnormal{ for all } I \neq \emptyset.    
\end{equation}
Due to the relations $l_i^2=0$ in $\CHQs(\GZip^{\mu})$,
the nonzero monomials of the Chern classes $l_i$ are precisely the $L_I$. We call the $L_I$ the Hodge-Chern monomials of $\GZip^{\mu}$. For $I \neq \emptyset$, let 
\begin{equation}
\label{eq-RI-rep}    
r_I=\bigoplus_{i \in I}r_i
\end{equation}
Then $\csf_{|I|}\Grif(G,\mu, r_I)=L_I$. Note that $r_I$ may not be a representation of $G$ but only of $\tilde G$, \eg if $G$ is adjoint. However $\Grif(G,\mu, r_I)$ is still well-defined thanks to~\S\ref{sec-tilde-G-reps}. Thus all the Hodge-Chern monomials $L_I$ are Hodge-Chern classes, \ie Chern classes of Griffiths-Hodge bundles. This uses in a crucial way that we consider Griffiths bundles relative $k$-representations $r$ rather than merely $\fp$-representations.
\end{subequations}
\subsection{Explicit strata-effectivity}
\label{sec-hilb-strat-eff}
\subsubsection{Zip strata}
\label{sec-hilb-EO-strata}
Given a subset $I \subset \ZZ/(d)$, let $I^c:=\ZZ/(d) \setminus I$ denote its complement. 
Identify the Weyl group $W=(\ZZ/2)^d$ with the power set  $\Pcal(\ZZ/(d))$ by associating to a subset $I \subset \ZZ/(d)$ the d-tuple with $1$ (resp. $0$) in coordinates $i \in I$ (resp. $j \in I^c$). Thus the zip strata of $\Xcal=\GZip^{\mu}$ are parameterized by $\Pcal(\ZZ/(d))$. 
The Zip stratum $\Xcal_I$ has dimension $|I|$. 
We write $\Xcal^I:=\Xcal_{I^c}$ and $\Xcal^i:=\Xcal^{\{i\}}$. 
 Wedhorn-Ziegler \cite{Wedhorn-Ziegler-tautological} show that every graded piece $ \CHQ^m(\GZip^{\mu})$ of the Chow ring $\CHQs(\GZip^{\mu})$ is spanned by classes of codimension $m$ strata closures. Hence $L_I$ admits an expression 
\begin{equation}
\label{eq-aij-linear-comb}
L_I=\sum_{J \subset \ZZ/(d), |J|=|I|} a(I,J) [\overline{\Xcal}^J] \quad \textnormal { with } a(I,J) \in \QQ.
\end{equation}
Our aim is to determine the coefficients $a(I,J)$ explicitly and to observe that they are all non-negative. 
\begin{theorem}
\label{th-eff-hilbert}
 For every $I \subset \ZZ/(d)$, the monomial $L_I \in \CHQ^{|I|}(\GZip^{\mu})$
is strata-effective.
\end{theorem}
More precisely, one reduces to the restriction of scalars case $\tilde G=\res^{\fpd}_{\fp}\SL(2)_{\fpd}$. Then:
\begin{theorem}
\label{th-precise-eff-hilbert} Assume that $\tilde G=\res^{\fpd}_{\fp}\SL(2)_{\fpd}$. Then:
\begin{enumerate}
\item
\label{item-th-aij-vanish}
An explicit characterization of when $a(I,J)=0$ in~\eqref{eq-aij-linear-comb} is given in~\ref{prop-aij-explicit}\ref{item-aij-vanish}.
\item
\label{item-th-aij-explicit}
When $a(I,J) \neq 0$, it admits the explicit formula~\ref{prop-aij-explicit}\ref{item-aij-explicit}. In particular, if $a(I,J) \neq 0$, then  $$a(I,J)=\frac{p^{e(I,J)}}{p^d+(-1)^{|I|}},
$$
where the exponent $e(I,J)$ satisfies $e(I,I)=d-|I|$ and $0 \leq e(I,J) \leq d-|I|-1$ for all $J \neq I$.
\end{enumerate}
\end{theorem}
\subsubsection{\texorpdfstring{$\textdbend$}{} Self-products} To prove~\ref{th-precise-eff-hilbert}, one may be tempted to simply multiply the codimension one formula~\ref{sec-codim-1} with itself $|I|$ times. Doing so is problematic because, while $[\overline{\Xcal}^i][\overline{\Xcal}^j]=[\overline{\Xcal}^{\{i,j\}}]$ for $i \neq j$ works out nicely,  one also encounters terms of the form $[\overline{X}^j]^2$ for various different $j$. These are most easily expressed as $-2l_jl_{j+1}$, which one must then express in terms of the strata classes. Moreover, by~\ref{th-eff-hilbert} $-2l_jl_{j+1}$ is strata anti-effective \ie has all coefficients non-positive when expressed in the basis of codimension two strata classes.
\subsubsection{Example: Codimension $1$} 
\label{sec-codim-1}
One has 
$$\left(p^d+(-1)^d\right)L_{\{i\}}:=\left(p^d+(-1)^d\right)l_i=\sum_{m \in \ZZ/(d)} p^{d-m-1}\overline{\Xcal}^{i+m}.
$$

\subsubsection{An example in codimension \texorpdfstring{$2$}{2}}
\label{sec-exa-24}
\addtocounter{equation}{-1}
\begin{subequations}
    
Consider the case where $d=5$ and $\tilde G=\res^{\fpd}_{\fp}\SL(2)_{\fpd}$.  Then: 
\begin{equation}
\label{eq-exa-24}
(p^5+1)L_{\{1,3\}}:=(p^5+1)l_1l_3=p^3[ \overline{\Xcal}^{\{1,3\}}]+p^2[\overline{\Xcal}^{\{2,3\}}]+p^2[\overline{\Xcal}^{\{1,4\}}]+p[\overline{\Xcal}^{\{2,4\}}]+p[\overline{\Xcal}^{\{0,1\}}]+[\overline{\Xcal}^{\{1,2\}}].
\end{equation}
\end{subequations}

\subsection{Application to Hilbert modular varieties}
\label{sec-hilbert-case}
\subsubsection{Hilbert modular varieties}
\label{sec-def-hilbert-modular} Assume that $(G,\mu)$, the special $k$-fiber $S_K$ and a smooth toroidal compactification of it $S_K^{\Sigma}$ arise as in~\ref{sec-intro-Hodge-type} from the Shimura datum of a Hilbert modular variety associated to a totally real number field $F$ of degree $d$ over $\QQ$. 
The canonical extension of the Hodge vector bundle $\Omega^{\can}$ admits a decomposition into line bundles
\begin{equation}
 \Omega^{\can}=\omega_{0}^{\can} \oplus \cdots \oplus \omega_{d-1}^{\can}   
\end{equation} corresponding to the $d$ embeddings of $F$ into $\RR$. 
Let $\iota : \tilde G \to G$ be the inclusion.
Let $e_i \in \Xsf^*(\tilde T)$ and $l_i \in \CHQ^1(\GZip^{\mu})$ as in~\S\ref{sec-hodge-chern-A1}. There are characters $\eta_i \in \Xsf^*(T)$ such that $\iota^*\eta_i=-e_i$. and that $\zeta^*\Lscr(\eta_i)=\omega_i^{\can}$. The first Chern class $\csf_1 \omega_i^{\can}=\zeta^{\Sigma,*}l_i$. 
\subsubsection{Ekedahl-Oort strata} The preimages $S_{K,I}:=\zeta^{-1}\Xcal_I$ and $S^{\Sigma}_{K,I}:=\zeta^{\Sigma,-1}(\Xcal_I)$ are the Ekedahl-Oort of $S_K$ and $S_K^{\Sigma}$ respectively. 

An immediate corollary of~\ref{th-eff-hilbert} and the smoothness of $\zeta^{\Sigma}$ is:
\begin{corollary}
    \label{cor-hilbert-var-strata-eff}
    For every $I \subset \ZZ/(d)$, the Hodge-Chern monomial $\zeta^{\Sigma,*}L_I = \prod_{i\in I}c_{1}(\omega_{i}^{\can})\in \Tsf^{|I|}(X^{\Sigma})$
is strata-effective on the Hilbert modular variety $S_K$. Explicitly, the coefficients of $\zeta^{\Sigma,*}L_I$ in terms of Ekedahl-Oort strata are given by~\ref{th-precise-eff-hilbert}. 
\end{corollary}
\begin{rmk} As is well-known, the case that $p$ is totally inert in $F$ corresponds to $G=\tilde G=\res^{\fpd}_{\fp}\SL(2)_{\fpd}$.
\end{rmk}

\subsection{Viewing the explicit strata-effectivity~\texorpdfstring{\ref{th-precise-eff-hilbert}}{} as a commutative algebra statement}
\label{sec-comm-alg}
\begin{subequations}
We explain how using the known, explicit description of the classes of strata in terms of $L_I$,~\ref{th-eff-hilbert} becomes a purely algebraic statement about the ring $R$~\eqref{eq-taut-hilbert}, which does not reference $\GZip^{\mu}$, Hilbert modular varieties or Chow groups. Our proof is then also purely algebraic, but see~\ref{sec-geometric} for some further geometric interpretation and motivation. 
\subsubsection{The $N_i$}
Assume that we are in the most interesting case where $\tilde G=\res^{\fpd}_{\fp}\SL(2)_{\fpd}$.
For all $i \in \ZZ/(d)$ and $I \subset \ZZ/(d)$ set 
\begin{equation}
N_i=pl_i-l_{i+1} \quad \textnormal{ and } \quad
N_I=\prod_{i \in I}N_i. 
\end{equation}
\subsubsection{Partial Hasse invariants}
Specialize the construction~\ref{sec-partial-Hasse} of partial Hasse invariants $\hsf^*f_{\lambda,w}$ on flag strata $\overline{\Ycal}_w$ to $G=\res^{\fpd}_{\fp}\SL(2)_{\fpd}$,  $w=w_0=\ZZ/d$ and $\lambda=e_i$. Note that $\GZip^{\mu}=\GF^{\mu}$ for $G$ of type $\Asf_1^d$. This gives partial Hasse invariants $h_i=\hsf^*f_{e_i,\ZZ/d} \in H^0(\GZip^{\mu}, \Lscr(pe_i-e_{i+1}))$. The zero schemes $\Zero(h_i)=\Xcal^i$ and $\Zero(\prod_{i \in I}h_i=\Xcal^I$. Hence  $N_J=[\overline{\Xcal}^J]$ in $\CHQs(\GZip^{\mu})$ for all $J \subset \ZZ/(d)$. Thus~\ref{th-precise-eff-hilbert} may be rephrased algebraically by replacing~\eqref{eq-aij-linear-comb} with 
\begin{equation}
\label{eq-LI-NJ}
L_I=\sum_{J \subset \ZZ/(d), |J|=|I|} a(I,J) N_J \quad \textnormal { with } a(I,J) \in \QQ,
\end{equation}
 forgetting (momentarily) the geometric interpretations of $R$, $L_I$ and $N_J$.
 \end{subequations}

\subsubsection{The matrix $A_m$ and its inverse}
The $m$th graded piece $\gr^m R$ of the $\QQ$-algebra $R$~\eqref{eq-taut-hilbert}  admits as basis all the $L_I$ with $|I|=m$.
Let $A_m \in \End( V)$ be the matrix of $N_J$ in the basis $(L_I)$.
That every $L_I$ is a linear combination of $N_J$ as in~\eqref{eq-LI-NJ} is equivalent to the invertibility of $A_m$. 
To see that $A_m$ is invertible without appealing to \cite{Wedhorn-Ziegler-tautological}, simply observe that $A_m$ has integer coefficients and its reduction mod $p$ is a permutation matrix, hence invertible. 
In this setting,~\ref{th-precise-eff-hilbert} explicitly describes the entries of the inverse $A_m^{-1}$. 
The claimed strata-effectivity amounts to the claim that all the entries of $A_m^{-1}$ are non-negative.    
\subsection{Reciprocity}
\label{sec-reciprocity}
\subsubsection{The positive counterparts of the $N_I$}
\label{def-PI}

For all $i \in \ZZ/(d)$ and all $I \subset \ZZ/(d)$ set 
$$ P_i:=pl_i+l_{i+1} \quad \textnormal{ and } \quad
P_I=\prod_{i \in I}P_i. $$
One has $N_I,P_I \in \Gr^{|I|}R$ and $N_IP_{I^c} \in \Gr^d R$.
\begin{rmk}[Difference of squares vanishing]
\label{rmk-diff-squares-vanish}
Note that
 $N_iP_i=p^2l_i^2-l_{i+1}^2=0$ in $R$ ~\eqref{eq-taut-hilbert}. It follows that $N_IP_J=0$ if $I \cap J \neq \emptyset$. 
 Together with~\ref{lem-NI-PJ}, this may be seen as the purely algebraic/combinatorial motivation for introducing the $P_I$ to study the $N_J$, without reverting to the algebro-geometric interpretation of the $L_I, N_J$ in terms of $\Tsf^*(X) \subset \CHQ^*(X)$. For some geometric motivation for the $P_I$, see~\ref{sec-PI-geometric}.
 \end{rmk}
\begin{lemma}
\label{lem-NI-PJ}
For all equinumerous subsets $I,J \subset \ZZ/(d)$, one has

$$ N_IP_{J^c}=\left\{ \begin{array}{cll} \left(p^d+(-1)^{|I|}\right) L_{\ZZ/(d)} & \textnormal{if} & I=J\\ 0 &\textnormal{if} & I\neq J\end{array} \right. \quad \textnormal{ in } \Gr^dR.$$      
\end{lemma}
\subsubsection{Proof}
If $I \neq J$, then $I \cap J^c \neq \emptyset$, so $N_IP_{J^c}=0$ by~\ref{rmk-diff-squares-vanish}.
Expanding $N_IP_{I^c}$ as a sum of $2^d$ monomials, the two terms which are nonzero constant multiples of $L_{\ZZ/(d)}=l_0 \cdots l_{d-1}$ are the ones with coefficient $p^d$ and $(-1)^{|I|}$. The other monomials vanish because each contains some $l_j$ with exponent $\geq 2$. \qed 
\subsubsection{\texorpdfstring{$\textdbend$}{!} Remark} Lemma~\ref{lem-NI-PJ} exhibits a duality between the $N_I$ and the $P_J$, which is further reinforced by the reciprocity~\ref{prop-reciprocity} below.  
\subsubsection{Decomposition into intervals}
\label{sec-intervals}
For all $a,b \in \ZZ/(d)$, let $[a,b]$ be the interval from $a$ to $b$ gotten by successively adding $1$. That is, $[a,b]:=\{a\}$ if $a=b$ and otherwise $[a,b]:=\{a,a+1, \ldots, b\}$. An interval $[a,b]$ is maximal in a subset $I \subset \ZZ/(d)$ if $a-1,b+1 \not \in I$. Every nonempty subset $I \subset \ZZ/(d)$ uniquely decomposes as a disjoint union of maximal intervals.

\begin{proposition}
 \label{prop-aij-explicit} Assume $\tilde G=\res^{\fpd}_{\fp}\SL(2)_{\fpd}$
The coefficient $a(I,J)$ in~\eqref{eq-aij-linear-comb} is given explicitly as follows: Let $J^c=\sqcup_t[a_t,b_t]$ be the decomposition of the complement $J^c$ into maximal intervals. Then:
\begin{enumerate}
\item
\label{item-aij-vanish} 
The coefficient $a(I,J) \neq 0$ if and only if $ |[a_t,b_t+1] \setminus I^c|=1$ for all $t$.
\item 
\label{item-aij-explicit}
If  $[a_t,b_t+1] \setminus I^c=\{a_t+e_t\}$ for all $t$,  set $e(I,J)=\sum_t e_t$. Then:
$$a(I,J)=\frac{p^{e(I,J)}}{p^d+(-1)^{|I|}}.
$$

 \end{enumerate}

\end{proposition}

\subsubsection{Coefficient reciprocity}
\label{sec-subsubsec-reciprocity}
\addtocounter{equation}{-1}
\begin{subequations} 
The main tool for proving~\ref{th-precise-eff-hilbert} is the following reciprocity satisfied by the coefficients $a(I,J)$. 
Express $P_I$ in the basis $(L_J)$: \begin{equation}
\label{eq-PI-LI}
P_I:=\sum_{J \subset \ZZ/(d), |J|=|I| } a^{\vee}(I,J) L_J 
\end{equation}
The coefficients $a(I,J)$ and $a^{\vee}(I,J)$ satisfy a reciprocity relative "transposed complementation" 
\begin{equation}
\label{eq-transp-complement}
(I,J)  \mapsto(J^c, I^c)    
\end{equation} and transposition of the negative $N$ terms with the positive $P$ terms: Up to the proportionality constant $p^d+(-1)^{|I|}$ which depends only on $(d, |I|)$, the coefficient $a(I,J)$ of $N_J$ in $L_I$ is given by the coefficient $a^{\vee}(J^c, I^c)$ of $L_{I^c}$ in $P_{J^c}$ \end{subequations}
 
\begin{proposition}[Reciprocity]
\label{prop-reciprocity}
For all equinumerous $I,J \subset \ZZ/(d),$ 

 $$\left( p^d+(-1)^d \right) a(I,J)=a^{\vee}(J^c,I^c).
 $$
\end{proposition}
\subsubsection{Remark}
\label{rmk-reduct-PI}
Since $P_I$ only involves positive signs, it is clear from their definition~\ref{eq-PI-LI} that the "dual" coefficients $a^{\vee}(I,J) \geq 0$. 
Hence~\ref{prop-reciprocity} implies~\ref{th-eff-hilbert} that $L_I$ is strata-effective. 
To prove the explicit description of $a(I,J)$ given in~\ref{th-precise-eff-hilbert} and~\ref{prop-aij-explicit}, we are reduced to explicitly computing $a^{\vee}(I,J)$, which we do in~\ref{lem-P-interval} and~\ref{cor-aij-vee-explicit}. 
\subsubsection{Proof of \texorpdfstring{~\ref{prop-reciprocity}}{}}
Let $K \subset \ZZ/(d)$ of the same cardinality as $I$. To prove~\ref{prop-reciprocity} for $(I,K)$, multiply both sides of~\ref{eq-aij-linear-comb} by $P_{K^c}$. If $J \neq K$, then $J \cap K^c \neq \emptyset$, so $N_JP_{K^c}=0$ by the difference of squares vanishing~\ref{rmk-diff-squares-vanish}. Hence the sum which forms the right-hand side of~\ref{eq-aij-linear-comb} simplifies to the term we want: $a_{I,K}N_KP_{K^c}$. By~\ref{lem-NI-PJ} this simplifies further to $a_{I,K}\left(p^d+(-1)^{|I|}\right)L_{\ZZ/(d)}$. 

The left-hand side of~\eqref{eq-aij-linear-comb} is similar, but easier: Expanding $P_{K^c}$, we get its expression as a linear combination of the monomials $L_{I'}$ of degree $|I'|=|I^c|=|K^c|=d-|I|$ with coefficients $a^{\vee}(K^c, I')$. If $I' \neq I^c$, then $I' \cap I \neq \emptyset$, so the product $L_IL_{I'}=0$ since $l_i^2=0$. Hence the left-hand side of~\eqref{eq-aij-linear-comb} simplifies to $a^{\vee}(K^c, I^c)L_{\ZZ/(d)}$.  
\
Comparing the left and right-hand sides, $a^{\vee}(K^c, I^c)L_{\ZZ/(d)}=a(I,K)L_{\ZZ/(d)}$. Since $L_{\ZZ/(d)} \neq 0$, the claimed reciprocity~\ref{prop-reciprocity} is true. 
\qed
\subsection{Explicit computation of the \texorpdfstring{$P_I$}{Ps} in terms of the \texorpdfstring{$L_J$}{Ls} }
\label{sec-PI-in-terms-of-LJ}
\begin{lemma}
\label{lem-P-interval}
For all $a,b \in \ZZ/(d)$,
    $$P_{[a,b]}= \sum_{r=0}^{|a-b|+1}p^rL_{[a,b+1]\setminus\{a+r\}}.$$
\end{lemma}

\subsubsection{Proof}
Expanding the product~\ref{def-PI} defining $P_{[a,b]}$ as a sum of $2^{|a-b|+1}$ monomials  and recalling that $l_i^2=0$ for all $i$ (so many of these monomials vanish), every nonzero   monomial $L_I$ which occurs is uniquely determined by the smallest index $r$ (if it exists) such that  $l_{a+r+1}$ is chosen  among $pl_{a+r},l_{a+r+1}$ from the $P_{a+r}=pl_{a+r}+l_{a+r+1}$ term.
Indeed then $pl_{a+s}$ is chosen  among $pl_{a+s},l_{a+s+1}$ for every $s<r$ by definition and $l_{a+s+1}$ is chosen  among $pl_{a+s},l_{a+s+1}$ for every $s>r$ by the nonvanishing assumption. This uniquely gives the monomial $L_{[a,b+1]\setminus\{a+r\}}$ with coefficient $p^r$. 
The case where no such $r$ exists corresponds uniquely to the monomial $L_{[a,b]}=L_{[a,b+1] \setminus\{b+1\}}$ with  coefficient $p^{|a-b|+1}$.
\qed
\begin{corollary}
\label{cor-aij-vee-explicit}
Let $I',J' \subset \ZZ/(d)$ be two nonempty equinumerous subsets. Let $J'=\sqcup_t[a_t,b_t]$ be the decomposition of $J'$ into maximal intervals~\eqref{sec-intervals}. 
The coefficient $a^{\vee}(J',I')$ in~\eqref{eq-PI-LI} is given as follows:
\begin{enumerate}
\label{item-aij-vee-vanish}
    \item The coefficient $a^{\vee}(J',I')\neq 0$ if and only if $ |[a_t,b_t+1] \setminus I'|=1$ for all $t$.
\item 
\label{item-aij-vee-explicit}
If  $[a_t,b_t+1] \setminus I'=\{a+e_t\}$ for all $t$, set $e^{\vee}(J',I')=\sum_t e_t$. Then: 
$$a^{\vee}(J',I')=p^{e^{\vee}(J',I')}.
$$

\end{enumerate}
    
\end{corollary}
\subsubsection{Remark} When $J'=J^c$ and $I'=I^c$, one has $e(I,J)=e^{\vee}(J^c,I^c)$.
\subsubsection{Proof of \texorpdfstring{~\ref{cor-aij-vee-explicit}}{}}
By the definition of maximal intervals, given $s \neq t$, the extended intervals $[a_s, b_s+1], [a_t,b_t+1]$ remain disjoint (in $\ZZ/(d)$; they no longer lie in $J'$). Hence
the decomposition of $P_{J'}$ according to intervals:
\begin{equation}
\label{eq-decomp-PJ-intervals}
P_{J'}=\prod_t P_{[a_t,b_t]}.   \end{equation}
So~\ref{cor-aij-vee-explicit} follows from~\ref{lem-P-interval}.
\qed
\subsubsection{Proof of~\ref{th-precise-eff-hilbert} and~\ref{prop-aij-explicit}}
As observed in~\ref{rmk-reduct-PI}, ~\ref{th-precise-eff-hilbert} and~\ref{prop-aij-explicit} follow by combining the results~\ref{lem-P-interval}~\ref{cor-aij-vee-explicit} on the $P_I$ and the dual coefficients $a^{\vee}(I,J)$ with the coefficient reciprocity~\ref{prop-reciprocity}. \qed
\subsubsection{Proof of~\ref{th-eff-hilbert}} By the results and constructions for equi-adjoint groups~\ref{prop-central-quotient}, \S\ref{sec-tilde-G-reps}, we reduce to the simply-connected case $G=\tilde G$. Let $G_n:=\res^{\FF_{p^{n}}}_{\fp}\SL(2)_{\FF_{p^{n}}}$. 
Then $G$ is a direct product $G=\prod_{j=1}^sG_{d_j}$ for some partition $d=\sum_{j=1}^s d_j$. If $\mu_j$ is the projection of $\mu$ onto $G_{d_j}$, then~\eqref{eq-product-GZipmu} gives the corresponding decomposition $$\GZip^{\mu}=\prod_{j=1}^s\GdjZip^{\mu_j}.$$For $I \subset \ZZ/d$, the Hodge-Chern monomial $L_I$ is a product of Hodge-Chern monomials $L_{I_j}$ for $\GdjZip^{\mu_j}$. Applying~\ref{th-precise-eff-hilbert} to each $L_{I_j}$ gives~\ref{th-eff-hilbert}, because the zip strata classes of $\GZip^{\mu}$ are the exterior products~\eqref{lem-Chow-GZip-prod} of those of the $\GdjZip^{\mu_j}$. \qed     

\subsection{Additional geometric information \texorpdfstring{\&}{} motivation}
\label{sec-geometric}
Each of the following statements~\ref{sec-geom-hodge-vector-bundle}-\ref{sec-PI-geometric} holds for all groups $G$ of type $\Asf_1^d$ on $\GZip^{\mu}$. To pull back  to a Hilbert modular variety $S_K$ via $\zeta$, specialize to the group $G$ which arises from the Shimura datum $\gx$. 
\subsubsection{The Hodge vector bundle}
\label{sec-geom-hodge-vector-bundle}
The Hodge vector bundle on $\GZip^{\mu}$ is $\Omega:=\fil^1\Wscr(r_{\ZZ/d}^{\vee})$. It has the same Chern classes as the Griffiths bundles $\Grif(G, \mu, r_{\ZZ/d}^{\vee})=\fil^1\Wscr(r_{\ZZ/d}^{\vee}) \oplus \Wscr(r_{\ZZ/d}^{\vee})$ and $\Grif(G, \mu, r_{\ZZ/d})$, since the Chern classes of $\Wscr(r_{\ZZ/d})$ and $\Wscr(r_{\ZZ/d}^{\vee})$ vanish.
When $G$ arises from $\gx$, its pullback $\zeta^{\sigma, *}\Omega=\Omega^{\can}$, the Hodge vector bundle on $S_K^{\Sigma}$.
\subsubsection{The top degree monomial \texorpdfstring{$L_{\ZZ/(d)}$}{}} The $l_i$ are also the Chern roots of the Hodge vector bundle $\Omega$ on $\GZip^{\mu}$. Hence the Hodge-Chern monomial $L_{\ZZ/(d)}=\csf_d(\Omega^{\can})$  is the top Chern class of $\Omega^{\can}$. When $G$ arises from the Shimura datum $\gx$, the pullbacks $\zeta^{\Sigma,*}l_i$ are the Chern rhoots of the Hodge vector bundle $\Omega^{\can}$ on $S_K^{\Sigma}$ and $\zeta^{\Sigma,*}L_{\ZZ/(d)}=\csf_d(\Omega^{\can})$.
\subsubsection{Geometric appearance of the \texorpdfstring{$P_I$}{PI}}
\label{sec-PI-geometric} Specializing the general construction of the partial Hasse invariants~\ref{sec-partial-Hasse} to a codimension $1$ stratum closure $\overline{\Xcal}^i=\overline{\Xcal}_{(\ZZ/d) \setminus\{i\}}$ gives the nowhere vanishing, partial Hasse invariant $$s_i:=\hsf^*f_{e_i,(\ZZ/d) \setminus\{i\}} \in H^0(\overline{\Xcal}^i,\Lscr(-pe_i-e_{i+1})).$$ 
Hence $[\overline{\Xcal}^i].P_i=0$ in $\CHQs(\GZip^{\mu})$. Pullback via $\zeta$ gives $[\overline{S}_K^{\Sigma, i}].P_i=0$ in the tautological ring $\Tsf^*(S_K^{\Sigma})$. This relation implies the relations $l_i^2=0$ which also follow from the Wedhorn-Ziegler isomorphism (\eqref{eq-Wedhorn-Ziegler-iso} in general,~\eqref{eq-wedhorn-ziegler-hilbert} in type $\Asf_1^d$).   The pullbacks $$\zeta^{\Sigma, *}s_i \in H^0(S_{K}^{\Sigma,i},(\omega_i^{\can})^{\otimes p}\omega_{i+1}^{\can})$$ are nowhere vanishing sections on codimension one EO strata which appear in the literature on Hilbert modular forms (mod $p$) \cf Emerton-Redduzi-Xiao \cite{Emerton-Reduzzi-Xiao} and generalize nowhere vanishing sections on the supersingular locus of modular curves studied by Serre, Edixhoven, Robert and others \cf \cite{Edixhoven-Serre-weight}.

\section{
Hodge-Chern classes II: Proportionality to powers of the Griffiths-Hodge line bundle}
\subsection{Powers of the Griffiths-Hodge line bundle} 
Recall the partial Hasse invariant cones $\Cscr_{\pha,w}^{P_0}$ associated to an intermediate parabolic $B \subset P_0 \subset G$~\ref{eq-def-partial-pha-cone}.
\begin{definition}
 Fix an integer $l$, $1 \leq l \leq |\Phi \setminus \Phi_{I_0}|/2$. A character $\lambda \in \Xsf^*(L_0)$ is a
\underline{length $l$ partial Hasse generator} for $\GF^{\mu,P_0}$ if $\lambda \in \Cscr^{P_0}_{\pha, w}$ for all $w \in \leftexp{I_0}{W}$ of length $l(w)=l$
 \end{definition}
In other words, $\lambda$ is a length $l$ partial Hasse generator for $\GF^{\mu, P_0}$ if $\Lscr(m\lambda)$ admits a partial Hasse invariant on every length $l$ stratum closure of $\GF^{\mu, P_0}$ of for some $m>0$. 
\begin{proposition}
\label{prop-length-hasse-eff} Fix an integer $l$, $1 \leq l \leq |\Phi \setminus \Phi_{I_0}|/2$. Assume that $\lambda \in \Xsf^*(L_0)$ is a length $l$ partial Hasse generator for $\GF^{\mu,P_0}$. Suppose that $\eta \in \Tsf_l(X^{P_0})$ is strata-effective. Then $\csf_1(\Lscr(\lambda)) \cdot \eta \in \Tsf_{l-1}(X^{P_0})$ is strata-effective.  
\end{proposition}
It may very well be that $\csf_1(\Lscr(\lambda)) \cdot \eta \in \Tsf_{l-1}(X^{P_0})=0$.
\subsubsection{Proof}
Write $\eta=\sum_{l(w)=l} a_w [\overline{X}^{P_0}_w]$. By assumption, for all $w$ of length $l$, the coefficient $a_w \geq 0$ and there exists a partial Hasse invariant $h_w \in H^0(\overline{X}^{P_0}, \Lscr(m_w \lambda))$ for some $m_w \geq 1$. 
Then $$m_w\csf_1(\Lscr(\lambda)) \cdot [X_w^{P_0}]=[\Zero(h_w)] $$   
is strata-effective. Hence 
$\csf_1(\Lscr(\lambda)) \cdot \eta$ is a non-negative linear combination of strata-effective classes. 
\qed
\begin{definition}
A character $\lambda \in \Xsf^*(L_0)$ is a
\underline{partial Hasse generator} for $\GF^{\mu,P_0}$ if $\lambda$ is a length $l$ partial Hasse generator for all $l$, $1 \leq l \leq |\Phi/ \Phi_{I_0}|/2$.  
\end{definition}
An induction on~\ref{prop-length-hasse-eff} gives:
\begin{corollary}
\label{cor-partial-hasse-gen-strata-effective}
Assume that $\lambda \in \Xsf^*(L_0)$ is a partial Hasse generator for $\GF^{\mu, P_0}$. Then $\csf_1(\Lscr(\lambda))^m$ is strata-effective for all $m$, $1 \leq m \leq |\Phi \setminus \Phi_{I_0}|/2$.
\end{corollary}
\begin{corollary}[\ref{th-intro-first-chern-griffiths-effective}]
\label{cor-grif-line-strata-eff}    
Assume that $r:G \to \GL(V)$ is an $\fp$-representation with central kernel. If $\mu$ is orbitally $p$-close~\ref{def-orb-p-close},  then $\csf_1(\Grif(\GZip^{\mu},r))^m$ is strata-effective for all $m \geq 0$.
\end{corollary}
\subsubsection{Proof}
By the main result of \cite{Goldring-griffiths-bundle-group-theoretic}, the Griffiths character $\grif(G,\mu,r)$ lies in the interior of the Griffiths-Schmid cone and it is orbitally $p$-close if and only if $\mu$ is. 
By the main theorem on group-theoretic Hasse invariants \cite[3.2.3]{Goldring-Koskivirta-Strata-Hasse}, the Griffiths-Hodge line bundle $\grif(\GZip^{\mu},r)$ is a partial Hasse generator for $\GZip^{\mu}$. \qed 
\subsection{Proportionality}
\label{sec-proportionality}
%
Part~\ref{item-grif-strat-eff} of the following is~\ref{th-grif-eff}.
\begin{theorem}
\label{th-proportionality}
Assume that $(\type(G), \type(L))$ is either $(\Xsf_n, \Xsf_{n-1})$ for some $\Xsf \in \{\Asf, \Bsf, \Csf, \Dsf\}$ or $(\Gsf_2, \Asf_{1})$. If $\type(G)=\Gsf_2$, assume $p \geq 5$.
Let
$r$ be a nontrivial, irreducible representation of $\tilde G$ of minimal dimension. Then: 
\begin{enumerate}
    \item 
\label{item-proportionality-grif}
Every Chern class of  $\Grif(\GZip^{\mu},r)$ is a strictly positive multiple  of $\csf_1 \Grif(\GZip^{\mu},r)$.
\item
\label{item-grif-strat-eff}
Every Chern class of  $\Grif(\GZip^{\mu},r)$ is strata-effective.
\end{enumerate}     
\end{theorem}
\begin{rmk}
\label{rmk-lowest-fundamental}
If $G$ in~\ref{th-proportionality} is not of type $\Asf_n$ ($n \geq 2)$ or $\Dsf_4$, then $\tilde G$ admits a unique nontrivial, irreducible representation of minimal dimension: the lowest dimensional fundamental representation of $\tilde G$. In the notation of Bourbaki's Planches \cite{bourbaki-lie-4-6}, it corresponds to the vertex $\alpha_1=e_1-e_2$ of the Dynkin diagram of $G$. 
If $G$ has type $\Asf_{n-1}$ and $n \geq 3$, then $\tilde G \cong \SL(n)$ has precisely two such representations: The inclusion $\SL(n) \hookrightarrow \GL(n)$ and its dual, which are the two fundamental representations corresponding to the pair of extremal vertices $e_1-e_2$ and $e_{n-1}-e_n$ swapped by the opposition involution $-w_0$. If $G$ has type $\Dsf_4$, it has three such representations: the standard/defining one and the two half-spin representations. They are the fundamental representations corresponding to the three extremal vertices $e_1-e_2, e_3-e_4, e_3+e_4$ permuted by triality.       
\end{rmk}
\begin{rmk}
Let $G=\Spin(2n)$ be the spin group of type $\Dsf_n$ for some $n \geq 2$. Assume that $\type(L)=\Dsf_{n-1}$. Let $r$ be the spin representation. Then not all Chern classes of $\Grif(G, \mu, r)$ are proportional to a power of its first Chern class. 
\end{rmk} 
\subsubsection{Proof of \texorpdfstring{~\ref{th-proportionality}}{}} By~\ref{cor-grif-line-strata-eff}, it suffices to show~\ref{item-proportionality-grif}. 
We explain the case where  $G$ has type $\Asf$. The other cases follow from similar computations. Assume that $\type(G)=\Asf_{n-1}$ (so re-index by $n-1$ in place of $n$). By~\S\ref{sec-funct-grif}, we may further assume that $G=\GL(n)$. In this case we may replace the Griffiths bundle with the Hodge vector bundle $\Omega$. It decomposes as $\Omega=\Omega_1 \oplus \Omega_{n-1}$ where $\Omega_i$ has rank $i$. 

Let $l_1, \ldots l_n$ be the Chern roots of the vector bundle $\Wscr(\id^{\vee})$ associated to the dual of the identity representation of $G=\GL(n)$. By \cite{Wedhorn-Ziegler-tautological}, the Chern classes of $\Wscr(\id^{\vee})$ vanish. One has $l_1=\csf_1(\Omega_1)$ and the Chern roots of $\Omega_2$ are $-l_2,\ldots ,-l_n$. Then~\ref{th-proportionality}\ref{item-proportionality-grif} in type $\Asf$ is given explicitly by~\ref{lem-proportionality-A}. \qed 

\begin{lemma}
\label{lem-proportionality-A}
For all $r \geq 1$, the Chern class $\csf_r(\Omega_{n-1})=l_1^r$.
\end{lemma}  
\subsubsection{Proof} By induction on $r$.
For $r=1$ this holds because $\sum_{i=1}^n l_i=\csf_1(\Wscr(\id^{\vee}))=0$. 

Then 
$$0=\sum_{1 \leq i_1 < \cdots <i_r \leq n} l_{i_1} \cdots l_{i_r}=
\sum_{2 \leq i_1 < \cdots <i_r \leq n} l_1l_{i_2} \cdots l_{i_r}+\sum_{2 \leq i_1 < \cdots <i_r \leq n} l_{i_1} \cdots l_{i_r}.
$$
$$=(-1)^{r-1}l_1\csf_{r-1}(\Omega_{n-1})+(-1)^r\csf_r(\Omega_{n-1})$$

By induction, $\csf_{r-1}(\Omega_{n-1})=l_1^{r-1}$. Hence $\csf_r(\Omega_{n-1})=l_1l_1^{r-1}=l_1^r$ as claimed. 
\qed
\section{ 
Effective tautological classes I: Divisors \texorpdfstring{\&}{} The cone conjecture}
\label{sec-cone-conj}
In this section we prove a more general version of \ref{th-pha-strata-eff} connecting strata effectivity of effective tautological classes with the cone conjecture. Using this we provide an explicit example of an effective tautological class which is not strata-effective for the Hilbert modular threefold at an inert prime.
\subsection{Theorem connecting strata-effectivity and the cone conjecture}

Recall the cones
\begin{equation*}
    \Cscr_{X,w}:=\{\lambda \in \Xsf^*(T) \ | \ H^0(\overline{Y}_{w}, \zeta_{Y}^*\Lscr(n\lambda)) \neq 0 \textnormal{ for some } n>0\}.
\end{equation*}
\begin{equation*}
    \Cscr_{\pha, w}:=\{\lambda \in \Xsf^*(T) \mid \Lscr(n\lambda)|_{\overline{\Ycal}_{w}} \text{ admits a partial Hasse invariant for some }n>0\}
\end{equation*} We prove the following for $w=w_0$:
\begin{theorem}[Theorem \ref{th-pha-strata-eff}]
    Let $(G,\mu)$ be a cocharacter pair. The following are equivalent.
    \begin{enumerate}
        \item For all $\zeta : X \rightarrow \GZip^\mu$ satisfying~\ref{conj-cone}\ref{item-smooth}-\ref{item-pseudo-complete}, $\Cscr_{X} = \Cscr_{\pha}$
        \item For all $\zeta : X \rightarrow \GZip^\mu$ satisfying~\ref{conj-cone}\ref{item-smooth}-\ref{item-pseudo-complete}, every effective class $\eta \in \Tsf^1(Y)$ on the flag space $Y/X$~\eqref{sec-intro-flag-space} $\eta$ is strata-effective. 
    \end{enumerate}
\end{theorem}
In order to prove this we need to introduce a new cone of injective sections. 
Define the cone
\begin{equation}
    \Cscr_{X,w,\inj} =\{\lambda \in \Xsf^*(T) \ | \textnormal{ there exists an injective section } s \in H^0(\overline{Y}_{w}, \zeta_{Y}^*\Lscr(n\lambda))  \textnormal{ for some } n>0\}.
\end{equation}
\begin{rmk}
    A priori, $\Cscr_{X,w,\text{inj}}$ could be strictly smaller than $\Cscr_{X, w}$ because there could be a non-zero section on $\overline{Y}_{w}$ which vanishes on some connected component, thus failing to be injective. One could hope to identify a property of $\zeta$ which implies $\Cscr_{X,w,\inj} = \Cscr_{X, w}$.
\end{rmk}
Since every partial Hasse invariant is injective, one has the inclusion of cones 
\begin{equation*}
    \Cscr_{\pha,w} \subset \Cscr_{X,w,\inj} \subset \Cscr_{X,w}.
\end{equation*} 
In order to prove Theorem \ref{th-pha-strata-eff} we prove the following statement with fixed $\zeta : X \rightarrow \GZip^\mu$. It is weaker than \ref{th-pha-strata-eff} in the sense that it restricts to the cone of injective sections, but more general in that it is a statement about all strata-closures.
\begin{theorem}
\label{th-pha-strata-eff-onstrata}
 Assume that $\zeta:X \to \GZip^{\mu}$ satisfies~\ref{conj-cone}\ref{item-smooth}-\ref{item-pseudo-complete}. Let $w \in W$. Then $\Cscr_{\pha,w}=\Cscr_{X,w,\inj}$ if and only if every effective Cartier divisor on the strata closure $\overline{Y}_{w}$ in the flag space $Y/X$~\eqref{sec-intro-flag-space} associated with an automorphic line bundle 
 is strata-effective.     
\end{theorem}
If $\overline{Y}_{w}$ is connected then the cones $\Cscr_{X,w} = \Cscr_{X,w,\inj}$ are equal and so \ref{th-pha-strata-eff-onstrata} gives a necessary and sufficient condition for $\Cscr_{\pha,w} = \Cscr_{X,w}$. Theorem \ref{th-pha-strata-eff} is a direct consequence of \ref{th-pha-strata-eff-onstrata}.
\begin{corollary}[Theorem \ref{th-pha-strata-eff}]
\label{cor-pha-strata-eff}
     Let $(G,\mu)$ be a cocharacter pair. The following are equivalent.
    \begin{enumerate}
        \item For all $\zeta : X \rightarrow \GZip^\mu$ satisfying~\ref{conj-cone}\ref{item-smooth}-\ref{item-pseudo-complete}, $\Cscr_{X} = \Cscr_{\pha}$
        \item For all $\zeta : X \rightarrow \GZip^\mu$ satisfying~\ref{conj-cone}\ref{item-smooth}-\ref{item-pseudo-complete}, every effective class $\eta \in \Tsf^1(Y)$ on the flag space $Y/X$~\eqref{sec-intro-flag-space} $\eta$ is strata-effective. 
    \end{enumerate}
\end{corollary}
\subsubsection{Proof}
    Apply \ref{th-pha-strata-eff-onstrata} with $w = w_{0}$ the longest element in the Weyl group. 
    
First, suppose that $\Cscr_{X} = \Cscr_{\pha}$ for all $\zeta$ satisfying~\ref{conj-cone}\ref{item-smooth}-\ref{item-pseudo-complete}. Suppose that $\zeta$ is one such morphism. By assumption, $\Cscr_{\pha} = \Cscr_{X,\inj} = \Cscr_{X}$  so every tautological effective Cartier divisor has strata-effective cycle class by \ref{th-pha-strata-eff-onstrata}. 
This is sufficient because every effective class $\eta \in \Tsf^1(Y)$ is the cycle class of an effective Cartier divisor since $Y$ is smooth. 

Conversely, suppose that for all $\zeta$ satisfying~\ref{conj-cone}\ref{item-smooth}-\ref{item-pseudo-complete}, every effective class $\eta \in \Tsf^1(Y)$ on the flag space $Y/X$~\eqref{sec-intro-flag-space} is strata-effective. Let $\zeta : X \rightarrow \GZip^\mu$ satisfy~\ref{conj-cone}\ref{item-smooth}-\ref{item-pseudo-complete}. 
For every connected component $X^{+}$ of $X$, $\zeta^{+} : X^{+} \rightarrow \GZip^{\mu}$ also satisfies~\ref{conj-cone}\ref{item-smooth}-\ref{item-pseudo-complete}. Hence every effective tautological divisor on $Y^{+}$ is strata-effective. 
Then \ref{th-pha-strata-eff-onstrata} implies that $\Cscr_{\pha} = \Cscr_{X^{+},\inj}$. On the other hand, $\Cscr_{X^{+},\inj}= \Cscr_{X^{+}}$  since every non-zero section on an integral scheme is injective. Note that $X^{+}$ is reduced because $\GZip^{\mu}$ is reduced and $\zeta$ smooth.
Since $\GZip^{\mu}$ and $\zeta$ are both smooth, so is $X$. In particular $X$ is normal. So the connected components of $X$ are irreducible. Hence $X^{+}$ is indeed integral. Pullback gives inclusions $\Cscr_{\pha} \subset \Cscr_{X} \subset \Cscr_{X^{+}}$ so $\Cscr_{\pha} = \Cscr_{X}$ as required. \qed
   

\subsubsection{Proof of~\ref{th-pha-strata-eff-onstrata}}
    Suppose first that $\Cscr_{\pha,w}=\Cscr_{X,\inj,w}$. Let $(\Lscr(\lambda),s)$ be an effective Cartier divisor on $\overline{Y}_{w}$. Then $s \in H^{0}(Y_{w},\zeta_{Y}^*\Lscr(\lambda))$ is injective. So by assumption $\zeta_{Y}^{*}\Lscr(n\lambda)$ admits a partial Hasse invariant on $\overline{Y}_{w}$ for some $n>0$. In particular, there is $\chi \in \CharT$ such that $\Lscr(n\lambda) = \hsf^{*}\Mscr_{w}(\chi)$. The cycle class of $(\zeta_{Y}^{*}\Lscr(\lambda)|_{\overline{Y}_{w}},s)$ in $\CHQ^*(Y)$ is $[\Zero(s)] = c_{1}(\zeta_{Y}^{*}\Lscr(\lambda)).[\overline{Y}_{w}]$. Now 
    \begin{equation*}
        \begin{aligned}
            n[\Zero(s)] &= c_{1}(\zeta_{Y}^{*}\Lscr(n\lambda)).[\overline{Y}_{w}] \\&= \zeta_{Y}^{*}\hsf^{*}(c_{1}(\Mscr_{w}(\chi)).[\overline{\Sbt}_{w}]) \\&= -\zeta_{Y}^{*}\hsf^{*}\sum_{\alpha \in E_{w}}\langle\chi,w\alpha^{\vee}\rangle[\overline{\Sbt}_{ws_\alpha}]
            \\&= -\sum_{\alpha \in E_{w}}\langle\chi,w\alpha^{\vee}\rangle[\overline{Y}_{ws_{\alpha}}]
        \end{aligned}
    \end{equation*}
    with the third equality holding by~\ref{le-cyc-cl-sbt}. Since $\zeta_{Y}^{*}\Lscr(n\lambda)$ admits a partial Hasse invariant on $\overline{Y}_{w}$ for some $n>0$, the coefficients $-\langle\chi,w\alpha^{\vee}\rangle$ are all non-negative by~\ref{le-admitpha}. Hence, $[\Zero(s)] = \frac{1}{n}c_{1}(\zeta_{Y}^{*}\Lscr(n\lambda)).[\overline{Y}_{w}]$ is strata-effective.
    
    
Conversely, suppose that every effective Cartier divisor on the strata closure $\overline{Y}_{w}$ associated with an automorphic line bundle is strata-effective. Let $s \in H^{0}(\overline{Y}_{w},\zeta_{Y}^{*}\Lscr(n\lambda))$ be an injective section. Then:
    \begin{equation*}
        \begin{aligned}
            \relax[\Zero(s)] &= c_{1}(\zeta_{Y}^{*}\Lscr(n\lambda)).[Y_{w}] \\&= \zeta_{Y}^{*}\hsf^{*}(c_{1}(\Mscr_{w}(\chi)).[\overline{\Sbt}_{w}]) \\&= \zeta_{Y}^{*}\sigma^{*}(-\sum_{\alpha \in E_{w}}\langle\chi,w\alpha^{\vee}\rangle[\overline{\Sbt}_{ws_{\alpha}}]) \\&= -\sum_{\alpha \in E_{w}}\langle\chi,w\alpha^{\vee}\rangle[\overline{Y}_{ws_{\alpha}}]
        \end{aligned}
    \end{equation*}
with the third equality holding by~\ref{le-cyc-cl-sbt}. By assumption $[\Zero(s)]$ is strata-effective so the coefficients $-\langle \chi,w\alpha^{\vee} \rangle$ are all non-negative. By~\ref{le-admitpha}, the pullback $\zeta_{Y}^{*}\Lscr(n\lambda)$ admits a partial Hasse invariant.
\qed
\subsection{An effective tautological class which is not strata-effective}
\label{sec-efftaut-notstrataeff}
In the case of a Hilbert modular variety $S_K$ at a prime $p$ inert in the totally real field $F$, the analogue of the Cone Conjecture~\ref{conj-cone} fails for codimension one EO strata. The cone of partial Hasse invariants $\Cscr_{\pha,w}$ is explicitly described in \cite[4.1.9]{Goldring-Koskivirta-global-sections-compositio}. On the other hand an explicit description of the cone of characters which give ample line bundles on the (open) Hilbert modular variety was conjectured by Liang-Xiao \cite{Tian-Xiao-quaternionic} and proven by Yang \cite{Yang-ampleness-hilbert}. These line bundles are ample on all positive codimension EO strata closures since these do not intersect the toroidal boundary. 
This provides us with a source of sections which are not, a priori, partial Hasse invariants.

Recall the notation~\ref{sec-intro-Hodge-type} for Hodge-type Shimura varieties and~\ref{sec-hilb-EO-strata} for the parameterization of strata for groups of type $\Asf_1^d$ by the power set $\Pcal(\ZZ/d)$. 
\begin{proposition}
\label{prop-hilbert-inert}
Let $d=3$. Let $S_K$ be a Hilbert modular threefold at a prime $p$ inert in $F$ and $S_K^{\Sigma}$ a smooth toroidal compactification. For all $I \subset \ZZ/(3)$ of size $2$, there exists $\lambda \in \CharT$ such that $\lambda \notin  \Cscr_{\pha,I}$, but $\lambda \in \Cscr_{S_K,I,\inj}$.
\end{proposition}
\subsubsection{Proof}
We describe the case $I=\{1,2\}$; the other two are handled similarly.
Let $\{e_{0},e_{1},e_{2}\}$ be the standard basis of $\QQ^{\ZZ/3}$. The second author and Koskivirta described $\Cscr_{\pha,J}$ for all $J \subset \ZZ/d$ in \cite[4.1.9]{Goldring-Koskivirta-global-sections-compositio}. In this case, with the choices made in \S\ref{sec-hilbert-case} we consider the $\QQ$-basis $\Bcal = \{\alpha_{0},\alpha_{1},\alpha_{2}\}$ of $\QQ^{\ZZ/3}$ with $\alpha_{0} = pe_{0} + e_{1}$, $\alpha_{1} = pe_{1} - e_{2}$ and $\alpha_{2} = pe_{1} - e_{0}$.
    $$\Cscr_{\pha,I} = \{m_{0}\alpha_{0} + m_{1}\alpha_{1} + m_{2}\alpha_{2} \mid m_{0} \in \QQ \text{, }m_{1},m_{2}\in \QQ_{\geq 0}\} \subset \QQ^{\ZZ/3}$$
    Given a vector $k_{0}e_{0}+k_{1}e_{1}+k_{2}e_{2} = m_{0}\alpha_{0} + m_{1}\alpha_{1} + m_{2}\alpha_{2} \in \QQ^{\ZZ/3}$ we have $$\begin{pmatrix}
        p & 0 & -1 \\ 1 & p & 0 \\ 0 & -1 & p
    \end{pmatrix}\begin{pmatrix}
        m_{0} \\ m_{1} \\ m_{2}
    \end{pmatrix} = \begin{pmatrix}
        k_{0} \\ k_{1} \\ k_{2}
    \end{pmatrix}$$
    Inverting this matrix gives:
    $$\frac{1}{p^{3}+1}\begin{pmatrix}
        p^{2} & 1 & p \\ -p & p^{2} & -1 \\ -1 & p & p^{2}
    \end{pmatrix}\begin{pmatrix}
        k_{0} \\ k_{1} \\ k_{2}
    \end{pmatrix} = \begin{pmatrix}
        m_{0} \\ m_{1} \\ m_{2}
    \end{pmatrix}.$$ So $k_{0}e_{0}+k_{1}e_{1}+k_{2}e_{2} \in \Cscr_{\pha,I}$ if and only if 
    \begin{enumerate}
        \item $(p^{3}+1)m_{1} = -pk_{0}+p^{2}k_{1}-k_{2} \geq 0$
        \item $(p^{3}+1)m_{2} = -k_{0} + pk_{1} + p_{2}k_{2} \geq 0$
    \end{enumerate}
Let $\lambda = k_{0}e_{0}+k_{1}e_{1}+k_{2}e_{2} \in \CharT$. By \cite[Theorem 1]{Yang-ampleness-hilbert}, the line bundle $\zeta^*\Lscr(\lambda)$ is ample on the open Hilbert modular variety $S_K \hookrightarrow S_K^{\Sigma}$ if and only if $pk_{0} > k_{2}$, $pk_{1} > k_{0}$ and $pk_{2} > k_{1}$ (note that under our conventions Frobenius acts as $i \mapsto i+1$). 
The EO strata closures  $\overline{S}_{K,J}$ of codimension $>0$ are complete. 
Since pullback under affine morphisms preserves ampleness, the ampleness of $\zeta^*\Lscr(\lambda)$  on $S_K$ implies the ampleness of $\zeta^{*}\Lscr(\lambda)|_{\overline{S}_{K,I}}$. The triple $(k_{0},k_{1},k_{2}) = (p^{3}-1,p^{2},p^{3}-1)$ satisfies $m_{1} < 0$ so $\lambda = k_{0}e_{0}+k_{1}e_{1}+k_{2}e_{2}$ does not lie in the partial Hasse invariant cone. 
However, $pk_{0} > k_{2}$, $pk_{1} > k_{0}$ and $pk_{2} > k_{1}$ so $\zeta^{*}\Lscr(\lambda)|_{\overline{S}_{K,I}}$ is ample. Hence some power $\zeta^{*}\Lscr(m\lambda)|_{\overline{S}_{K,I}}$ is globally generated and, in particular, admits an injective global section. Hence, $\lambda \in \Cscr_{S_K,I, \inj}$.
\qed

Combining~\ref{prop-hilbert-inert} with~\ref{th-pha-strata-eff-onstrata} gives:
\begin{corollary}
\label{cor-tauteff-notstrataeff}
In a Hilbert modular threefold $S_K^{\Sigma}$, there exist tautological curves  which are not strata-effective.
\end{corollary}
\begin{rmk}
 The non-strata-effective tautolgical curves in~\ref{cor-tauteff-notstrataeff} are contained in a codimensione one stratum of $S_K^{\Sigma}$. Assume that $\type(G)=\Asf_1^3$ and that $\zeta:X \to \GZip^{\mu}$ satisfies~\ref{conj-cone}\ref{item-smooth}-\ref{item-pseudo-complete}. 
 Combining the proof \cite[4.2.4]{Goldring-Koskivirta-global-sections-compositio} of the Cone Conjecture for groups of type $\Asf_1^d$ with the strata-effectiving of generically-ordinary curves~\ref{th-gen-ord-hilb} shows: The class of a closed subscheme $Z \subset X$ is strata effective unless $Z$ admits a one-dimensional irreducible component which densely intersects a codimension one stratum.    
\end{rmk}
\begin{rmk} In \cite{Brunebarbe-Goldring-Koskivirta-Stroh-ampleness} it is shown that $\Lscr(\lambda)$ is ample modulo the boundary on the flag space of a Hodge-type Shimura variety if $\lambda$ lies in the interior of the Griffiths-Schmid cone and $\lambda$ is orbitally $p$-close~\ref{def-orb-p-close}. The weight $(p^3-1, p^2, p^3-1)$ lies in the Griffiths-Schmid cone but it is not orbitally $p$-close. 
\begin{rmk}
 Yang's proof of the ampleness result recalled above \cite{Yang-ampleness-hilbert}, hence also~\ref{prop-hilbert-inert}, uses that $X$ is a Hilbert modular variety rather than only assuming that there is a smooth surjective Zip period map $X \to \GZip^{\mu}$. It remains to be seen whether~\ref{prop-hilbert-inert} holds for all $X$ admitting a smooth, surjective Zip period map.    
\end{rmk}    
\end{rmk}


\section{Effective tautological classes II: curves}
\label{sec-curve-strata-eff}

Let $\zeta: X\to \GZip^{\mu}$ be a zip period map. 
Let $\zeta_Y: Y \to \GF^{\mu}$ be the associated flag morphism~\eqref{sec-strat-zip-period}. 
\subsection{Generically \texorpdfstring{$w$}{w}-ordinary tautological curves}
\subsubsection{Tautological curves}
A tautological curve $C$ in $X$ (resp. $Y$) is a curve $C \subset X$ (resp. $Y$) whose class $[C] \in \CH_{1,\QQ}(X)$ (resp. $[C] \in \CH_{1,\QQ}(Y)$) is tautological, \ie lies in  $\Tsf_1(X)=\Tsf^{\dim X-1}$ (resp. $\Tsf_1(Y)=\Tsf^{\dim Y-1}$). 
\subsubsection{Generically $w$-ordinary curves}
\label{sec-generically-ordinary}
For $w \in W$, a curve $C \subset X$ (resp. $C \subset Y)$ is \underline{generically $w$-ordinary} if its intersection with the $w$-stratum\footnote{This definition makes sense without any assumption on $\zeta$, in which case the $X_w$ (resp. $Y_w$) may not satisfy the stratification property.} $X_{w}$ (resp. $Y_w$) is open, dense in $C$. For $w=w_0$, we say $C$ is \underline{generically ordinary}. 
\subsubsection{$\textdbend$}
Let $X=S_K$ be the special fiber of a Shimura variety~\ref{sec-intro-Hodge-type} and let $C \subset Y$ be a curve in its flag space which is generically ordinary in the sense just defined~\eqref{sec-generically-ordinary}. The image of $C$ in $S_K$ may not be contained in the ordinary locus. 
\subsubsection{Effective and strata-effective $w$-ordinary cones} 
Let $\Eff^1_{\taut,\word}(Y)$ be the cone of effective, tautological and generically $w$-ordinary curves.
Let $\Eff^1_{\strat,\word}(Y)$ be the sub-cone of strata-effective curves in $\Eff^1_{\taut,\word}(Y)$.

\subsection{Dual partial Hasse cone criterion}

\subsubsection{Dual cones}
Recall that $\Ycal:=\GF^{\mu}$ and that $\CH_{\QQ,0}(\Ycal)$ is one-dimensional, spanned by the class $[\Ycal_1]$ of the minimal stratum $\Ycal_1$, which is closed. 
Let $\Cscr \subset \Tsf^i(Y)$ be a cone. Define its dual $\Cscr^{\vee} \subset \Tsf_{i}(Y)$ by
$$\Cscr^{\vee}=\{\eta \in \Tsf_{i}(Y) | \eta \cdot \theta \textnormal{ is a non-negative multiple of } [Y_1] \textnormal{ for all } \theta \in \Tsf^i(Y) \}.
$$ 
Let $\csf_1 \Cscr_{\pha}(Y)=\{\csf_1\Lscr(\lambda) | \lambda \in \Cscr_{\pha}\} \subset \Tsf^1(Y)$.
\begin{proposition} 
\label{prop-dual-zip}
\addtocounter{equation}{-1}
\begin{subequations}
Assume that $X$ is a proper $k$-scheme and that $\zeta$ satisfies~\ref{conj-cone}\ref{item-smooth}. Suppose that 
\begin{equation}
 \label{eq-pha-dual-eff-cone}   
 \csf_1\Cscr_{\pha}^{\vee}(Y) \subset \Eff^1_{\strat,\word}(Y).
\end{equation}
Then every tautological, generically $w$-ordinary curve $C \subset Y$ is strata-effective:
\begin{equation}
\label{eq-eff-taut-strat-cone}
 \Eff^1_{\taut,\word}(Y) \subset \Eff^1_{\strat,\word}(Y).  
\end{equation}
\end{subequations}
\end{proposition}
\subsubsection{Proof}
We claim that $$\Eff^1_{\taut,\word}(Y) \subset \csf_1\Cscr_{\pha}^{\vee}(Y). $$
Let $\lambda \in \Cscr_{\pha}$. By definition, there exists $m>0$ and a partial Hasse invariant  $f \in H^0(\overline{\Ycal}_w, \Lscr(m\lambda))$. The intersection product $\csf_1\Lscr(m\lambda).C=m\csf_1\Lscr(\lambda).C$ on $Y$ is the class of the zero scheme $\Zero(\zeta_Y^*f|_{C})$ of a  restricted to $C$, hence is non-negative. 
\qed

\begin{rmk}
 The interest of~\ref{prop-dual-zip} is that its cone inclusion hypothesis is intrinsic to $\Xcal=\GZip^{\mu}$ and does not depend on $(X,\zeta)$ in the following sense: By the generalized Wedhorn-Zieglier isomorphism~\ref{eq-wedhorn-ziegler-generalized}, the strata-effective and dual partial Hasse cones $\Eff^1_{\strat,\word}(Y), \csf_1\Cscr_{\pha}(Y)$ are independent of $(X,\zeta)$ in the sense that $\csf_1\Cscr_{\pha}(Y)=\zeta^*\csf_1\Cscr_{\pha}(\Ycal)$ and $\Eff^1_{\strat,\word}(Y)=\zeta^*\Eff^1_{\strat,\word}(\Ycal)$ for all $\zeta$ as in~\ref{prop-dual-zip}.
 The validity of~\eqref{eq-pha-dual-eff-cone} for a triple $(G,\mu, w)$ is an explicit computation in terms of Chevalley's formula.     
\end{rmk}

\subsection{Type \texorpdfstring{$\Asf_1^d$}{a power of A1} and Hilbert modular varieties}
The following is~\ref{th-intro-curves} for type $\Asf_1^d$.
\begin{theorem}
\label{th-gen-ord-hilb}
Let $(G,\mu)$ be as in \ref{sec-def-hilbert-modular}. Suppose that $X$ is projective over $k$ and $\zeta : X \rightarrow \GZip^{\mu}$ is smooth and surjective. Then every tautological and generically ordinary curve $C \subset X$ is strata-effective.
\end{theorem}
In fact, we only need to assume that the zero-dimensional stratum $X_{e}$ has non-zero cycle class in $\Tsf_{0}(X)$. This is immediate if $X$ is a proper scheme.
\subsubsection{Proof}
Suppose that $C \subset X$ is generically ordinary. Since $C$ is generically ordinary, the partial Hasse invariants which cut out the codimension $1$ EO strata are injective when restricted to $C$. 
    \begin{equation}
        \label{eq-gen-ord-hilb}[C].(pl_{i}-l_{i+1}) \geq 0
    \end{equation}
     for all $i = 1,\ldots,d$. Since $C$ is tautological $[C] = \sum_{i=1}^{d}a_{i}L^{i}$ where $L^{i} := \prod_{j \in \mathbb{Z}/(d)\setminus\{i\}}l_{j}$. It suffices to show that $a_{i} \geq 0$ by \ref{th-eff-hilbert}. For $i = 1,\ldots,d$, \ref{eq-gen-ord-hilb} implies that $[C].(pl_{i}-l_{i+1}) = (pa_{i}-a_{i+1})l_{1}\ldots l_{d} \geq 0$. This immediately implies that $a_{i} \geq 0$ for all $i=1,\ldots,d$.
\qed
\subsection{The split \texorpdfstring{$\Asf_1^d$}{power of A1}-case and Hilbert modular varieties at a split prime}
\label{sec-split-hilbert}
Using the same method one can show the following:
\begin{theorem}
\label{th-hilbert-split-eff-taut}
Assume that $\zeta:X \to \GZip^{\mu}$ is a smooth and surjective. Suppose that $G^{\ad} \cong \PGL(2)_{\fp}^d$ for some $d \geq 1$ and that the centralizer $L:=\cent(\mu)$ is a torus. Then every effective tautological cycle $\eta \in \Tsf^*(X)$ is strata-effective. 
\end{theorem}
\begin{corollary}
\label{cor-hilbert-split-eff-taut}
Let $S_K$~\eqref{sec-intro-Hodge-type} be the mod $p$ special fiber of a Hilbert modular variety associated to a totally real field $F$~\eqref{sec-def-hilbert-modular}. Assume $p$ splits completely in $F$. Every effective tautological cycle $\eta \in \Tsf^*(S_K)$ is strata-effective.
\end{corollary}

\subsection{Type \texorpdfstring{$\Csf_2$}{C2}}
The following is~\ref{th-intro-curves} for type $\Csf_2$.
\begin{theorem}
\label{th-gen-ord-Sp4}
Assume that $\type(G)=\Csf_2$ and that $\mu^{\ad} \in \Xsf_*(G^{\ad})$ is minuscule. Suppose that $X$ is projective over $k$ and that  $\zeta:X \to \GZip^{\mu}$ is smooth and surjective. Then every generically ordinary, tautological curve $C \subset Y$ is strata-effective.   
\end{theorem}
\begin{lemma}
\label{lem-taut-Sp(4)}
Under the hypotheses of~\ref{th-gen-ord-Sp4}, 
$\Tsf^*(X)=\QQ[\lambda_1, \lambda_2]/(\lambda_1^2-2\lambda_2)$
\end{lemma}
\subsubsection{Proof}
This is seen in at least two ways: It is a special case of the Brokemper and Wedhorn-Ziegler isomorphisms ~\eqref{sec-intro-wedhorn-ziegler}~\eqref{eq-brokember-iso},~\eqref{eq-iso-L-conj-cpt-dual}
coupled with the well-known formula for the cohomology of the compact dual $\XX^{\vee}$. It also follows from relations in~\ref{fig-C2-min-taut}.   
\subsubsection{Proof of \texorpdfstring{~\ref{th-gen-ord-Sp4}}{}}
Due to the relation $l_1^2+l_2^2=0$ in $\Tsf^2(Y)$, the class $[C]=al_1^2l_2+bl_1l_2^2$ for some $a,b \in \QQ$. By~\ref{fig-C2-min-taut}, the line bundles $\Lscr(-(p-1),-(p-1))$ and $\Lscr(1,-p)$ admit partial Hasse invariants on $Y$. Hence  $(p-1)(l_1+l_2)[C] \geq 0$ and  $(pl_2-l_1)[C] \geq 0$. Expanding and using $l_1^2+l_2^2=0$ again gives $(bp+a)l_1l_2^3 \geq 0$ and $(b-a)l_1l_2^3 \geq 0$. By~\ref{fig-C2-min-taut}, $[Y_1]=(p^4-1)l_1l_2^3$. Hence $l_1l_2^3 \geq 0$. So $bp+a \geq 0$ and $b-a \geq 0$.  

On the other hand, using the formulas~\ref{fig-C2-min-taut},  $[\overline{Y}_{(12)}]$ (resp. $[\overline{Y}_{\sgn_2}]$) is a positive multiple of $l_1^2l_2+l_1l_2^2$ (resp. $-pl_1^2l_2-l_1l_2^2$). Inverting the matrix $\left(\begin{array}{cc}
    -1 & 1 \\
     -p& 1
\end{array} \right)$, a class $al_1^2l_2+bl_1l_2^2$ is a non-negative $\QQ$-linear combination of the strata classes $[\overline{Y}_{(12)}]$, $[\overline{Y}_{\sgn_2}]$ if and only if $b-a \geq 0$ and $bp-a \geq 0$. The curve $C$ is strata-effective because $bp+a \geq 0$ and $b-a \geq 0$ implies that $b-a \geq 0$ and $bp-a \geq 0$.   
\qed
\begin{figure}[ht]     
\caption{Type $\Csf_2$: Partial Hasse invariant-induced relations between strata classes}

\label{fig-C2-min-taut}
 \centerline{
$\xymatrixcolsep{-3pc}\xymatrixrowsep{2pc}\xymatrix{  &
\framebox{$\begin{array}{c}w_0=-1  \\ \leftrightarrow (14)(23) \end{array} : 1$ }  
\ar@{-}[rd]^(.55){-l_1+pl_2} 
\ar@{-}[ld] 
\ar@<-3pt>@{}[ld]_(.55){(p-1)(l_1+l_2)}  
&
\\ 
\framebox{$\begin{array}{c}
s_{\alpha}s_{\beta}s_{\alpha}=\sgn_1  \\ \leftrightarrow (14)
\end{array} : \begin{array}{l}
     (p-1)(l_1+l_2)  \\
   =(p-1)\lambda_1   
\end{array} $ } 
\ar@{-}[rrd]
\ar@<5pt>@{}[rrd]^(.4){\frac{1}{2}((p-1)l_1+(p+1)l_2)}
\ar@{-}[d]_-{\frac{1}{2}(-(p+1)l_1+(p-1)l_2)}
& 
&
\framebox{$\begin{array}{c}
 s_{\beta}s_{\alpha}s_{\beta}=-(12)      \\
   \leftrightarrow   (13)(24)
\end{array} : -l_1+pl_2$ } 
\ar@{-}[lld]_(.3){(p-1)l_2}
\ar@{-}[d]^-{(p-1)l_1}
\\
\framebox{$\begin{array}{c}
    s_{\beta}s_{\alpha}= \\
    \sgn_2 \circ (12)   \\
     \leftrightarrow (1243)
\end{array} : \begin{array}{l}
(p-1)(pl_2^2-l_1l_2)=       \\
(p-1)(-l_1l_2+\frac{1}{2}((p-1)l_2^2-(p+1)l_1^2))  \\ \hline   [\Rightarrow l_1^2+l_2^2=0] 
\end{array}  $} 
\ar@{-}[rrd]^(.65){(p+1)l_2}
\ar@{-}[d]_-{(p-1)l_1-(p+1)l_2}
&
&
\framebox{$\begin{array}{c}
     s_{\alpha}s_{\beta}= \\ \sgn_1 \circ (12)  \\
     \leftrightarrow (1342) 
\end{array} :  \begin{array}{l} (p-1)(pl_1l_2-(p-1)l_1^2)= \\ 
(p-1)(pl_1l_2+\frac{1}{2}((p-1)l_1^2+(p^2-1)l_2^2 ))
\\ \hline   [\Rightarrow l_1^2+l_2^2=0]
\end{array}$ }
\ar@{-}[lld]
\ar@<-3pt>@{}[lld]_(.7){(p+1)l_1+(p-1)l_2}
\ar@{-}[d]^-{-(p+1)l_1}
\\  
\framebox{$\begin{array}{c}
   s_{\alpha}=(12)    \\
   \leftrightarrow (12)(34)  
\end{array} :  
\begin{array}{l} (p-1)(p^2+1)(l_1^2l_2+l_1l_2^2) \\
=(p-1)(p^2+1)\lambda_1\lambda_2
\end{array}$
}
\ar@{-}[rd]_-{(p+1)l_2}
\ar@{-}[d]_-{(p+1)(l_1+l_2)}
&
&
\framebox{$\begin{array}{c}
    s_{\beta}=\sgn_2   \\
      \leftrightarrow (23)
\end{array} :  \begin{array}{l} -(p^2-1)(pl_1^2l_2+l_1l_2^2)
\end{array}$
}
\ar@{-}[ld]^-{pl_1-l_2}
\ar@{-}[d]^-{-l_1-pl_2}
\\
\framebox{$\emptyset :  \begin{array}{l} (p^4-1)(l_1^3l_2+2l_1^2l_2^2+l_1l_2^3) \\ \hline [\Rightarrow l_1^2l_2^2=0]
\end{array}$
}
& 
\framebox{$1 :\begin{array}{l}
     (p^4-1)(l_1^2l_2^2+l_1l_2^2)  \\
      =(p^4-1)l_1l_2^3 \\ \hline   [\Rightarrow l_1^2l_2^2=0]
\end{array}$}
&
\framebox{$\emptyset :  \begin{array}{l} (p^4-1)l_1^2l_2^2=0 \\ \hline [\Rightarrow l_1^2l_2^2=0]
\end{array}$
}
} $}
\end{figure}
\subsection{Type \texorpdfstring{$\Asf_2$}{A2}-unitary}
$\Omega = \Vcal(\std_{2}\oplus \std_{1}^{\vee})$ has weights $\{(-1,0,0),(0,-1,0),(0,0,1)\}$.
\begin{lemma}
    $\Tsf^*(Y) = \QQ[l_1,l_2]/(l_1^2 + l_2^2 + l_{1}l_2, l_{1}l_2^2 + l_1^{2}l_2)$
\end{lemma}
\subsubsection{Proof}
We have the relation $1 = (1-(l_1+l_2)t+l_{1}l_{2}t^2)(1+l_{3}t)$ since $\std_{2}\oplus \std_1 = \std$ is pulled back from a representation of $G$ so the chern classes of its automorphic vector bundle vanish.
\qed

These relations are witnessed concretely by partial Hasse invariants, see diagram \ref{fig-A2-unitary-min-taut}.
\begin{figure}[ht]     
\caption{Type $\Asf_2$-unitary: Partial Hasse invariant-induced relations between strata classes}

\label{fig-A2-unitary-min-taut}
 \centerline{
$\xymatrixcolsep{-3pc}\xymatrixrowsep{2pc}\xymatrix{  &
\framebox{$\begin{array}{c}w_0 = (13) \end{array} : 1$ }  
\ar@{-}[rd]^(.55){(p-1)l_1+pl_2} 
\ar@{-}[ld] 
\ar@<-3pt>@{}[ld]_(.55){(p-1)l_2-l_1}  
&
\\ 
\framebox{$\begin{array}{c}
  (123)
\end{array} : \begin{array}{l}
     (p-1)l_2-l_1
\end{array} $ } 
\ar@{-}[rrd]
\ar@<5pt>@{}[rrd]^(.4){(p+1)l_2}
\ar@{-}[d]_-{(p-1)l_1-l_2}
& 
&
\framebox{$\begin{array}{c}
       (132)
\end{array} : (p-1)l_1+pl_2$ } 
\ar@{-}[lld]_(.3){pl_1+(p-1)l_2}
\ar@{-}[d]^-{-(p+1)l_1}
\\
\framebox{$\begin{array}{c}
      (12)
\end{array} : \begin{array}{l}
((p^2-1)+1)l_{1}l_2 - (p-1)(l_1^2+l_2^2)=       \\
(p^2+(p-1)^2)l_{1}l_2 + p(p-1)(l_1^2+l_2^2)= \\
(p^2-p+1)l_{1}l_2\\ \hline   [\Rightarrow l_1^2+l_2^2+l_{1}l_2=0] 
\end{array}  $} 
\ar@{-}[rd]^(.65){-(p+1)l_1}
&
&
\framebox{$\begin{array}{c}
      (23)
\end{array} :  \begin{array}{l} -(p+1)(pl_{1}l_2 + (p-1)l_1^2)= \\ 
(p+1)(-l_{1}l_2+(p-1)l_2^2)
\\ \hline   [\Rightarrow l_1^2+l_2^2+l_{1}l_2=0]
\end{array}$ }
\ar@{-}[ld]
\ar@<-3pt>@{}[ld]_(.7){(p-1)l_1-l_2}
\\
& 
\framebox{$1 :\begin{array}{l}
     -(p+1)(p^2-p+1)l_1^2l_2  \\
      = -(p+1)((p-1)^{2}l_1^2l_2 - pl_{1}l_2^2)\\ \hline   [\Rightarrow l_1^2l_2+l_{1}l_2^2 = 0]
\end{array}$}}
$}
\end{figure}

\begin{theorem}
\label{th-gen-ord-unip21-inert}
Let $G$ be non-split of type $\Asf_2$. Suppose that $X$ is projective over $k$ and that  $\zeta:X \to \GZip^{\mu}$ is smooth and surjective. Then every generically ordinary, tautological curve $C \subset Y$ is strata-effective. 
\end{theorem}
\subsubsection{Proof}
    We make use of \ref{fig-A2-unitary-min-taut}. First we claim that a class $\alpha = al_{1}l_2 + bl_2^2 \in \Tsf_1(Y)$ is strata-effective if and only if it satisfies $(p-1)a+b \geq 0$ and $b \geq 0$. Writing $\alpha = A[\overline{Y}_{(12)}]+B[\overline{Y}_{(23)}]$, it is strata-effective if and only if $A\geq 0$ and $B\geq 0$. From the diagram \ref{fig-A2-unitary-min-taut},
    \begin{equation*}
        \begin{pmatrix}
            [\overline{Y}_{(12)}] \\ [\overline{Y}_{(23)}]\end{pmatrix} = 
            \underbrace{\begin{pmatrix}
                p^2-p+1 & 0 \\ -(p+1) & (p+1)(p-1)
            \end{pmatrix}}_{=:M}\begin{pmatrix}
                l_{1}l_2 \\ l_2^2
            \end{pmatrix}
    \end{equation*} By linear algebra, 
    \begin{equation*}
        \begin{pmatrix}
            A \\ B
        \end{pmatrix} = \underbrace{\frac{1}{\det M}\begin{pmatrix}
           (p+1)(p-1)  & p+1 \\ 0 & p^2-p+1
        \end{pmatrix}}_{(M^T)^{-1}}\begin{pmatrix}
            a \\ b
        \end{pmatrix}
    \end{equation*}
    Whence, the claim.
    Now suppose that $C\subset Y$ is generically ordinary and tautological. Any partial Hasse invariant cutting out a codimension $1$ EO stratum-closure $\overline{Y}_{w}$ will be injective when restricted to $C$ since it does not vanish on any component of $C$ (since $C$ generically ordinary). This gives that $((p-1)l_1+pl_2)[C]$ and $((p-1)l_2-l_1)[C]$ are effective. Writing $[C] = al_{1}l_2 + bl_2^2$ obtains $a+(p-1)b \geq 0$ and $pa-b \geq 0$. Note that this doesn't imply $[C]$ strata-effective since $(p,-1)$ satisfies these inequalities. To complete the proof we require additional positivity from somewhere. In this case we utilise the fact that the Hodge line bundle $\omega/X$ is nef. Arbitrary pullback of nef bundles is nef so $\pi^{*}\omega/Y$ is also nef. By definition of nef $c_{1}(\pi^{*}\omega).[C] = 2(l_1+l_2)(al_{1}l_2+bl_2^2)\geq 0$ for all closed curves $C \subset Y$. This yields the additional inequality $b \geq 0$ as required to finish the proof that $[C]$ is strata-effective.
\qed
\subsection{Type \texorpdfstring{$\Asf_2$}{A2}-split}
$\Omega = \Vcal(\std_{2}\oplus \std_{1}^{\vee})$ has weights $\{(-1,0,0),(0,-1,0),(0,0,1)\}$.
\begin{lemma}
    $\Tsf^*(Y) = \QQ[l_1,l_2]/(l_{1}^2+l_{2}^2+l_{1}l_2,l_1^{2}l_2+l_{1}l_2^2)$ 
    
\end{lemma}
We have the relations $l_1^2+l_2^2+l_{1}l_2 = 0$ and $l_1^{2}l_2+l_{1}l_2^2 = 0$ since $\std_{2}\oplus \std_1 = \std$ is pulled back from a representation of $G$ so the chern classes of its automorphic vector bundle vanish.

\begin{figure}[ht]     
\caption{$\Asf_{2}$-split: Partial Hasse invariant-induced relations between strata classes}
\label{fig-A2-slit-min-taut}
 \centerline{
$\xymatrixcolsep{-3pc}\xymatrixrowsep{2pc}\xymatrix{  &
\framebox{$\begin{array}{c}w_0 = (13) \end{array} : 1$ }  
\ar@{-}[rd]^(.55){pl_2-l_1} 
\ar@{-}[ld] 
\ar@<-3pt>@{}[ld]_(.55){(p-1)l_1+(p-1)l_2}  
&
\\ 
\framebox{$\begin{array}{c}
  (123)
\end{array} : \begin{array}{l}
     (p-1)l_1+(p-1)l_2
\end{array} $ } 
\ar@{-}[rrd]
\ar@<5pt>@{}[rrd]^(.4){pl_1+(p+1)l_2}
\ar@{-}[d]_-{-(p+1)l_1-l_2}
& 
&
\framebox{$\begin{array}{c}
       (132)
\end{array} : pl_2-l_1$ } 
\ar@{-}[lld]_(.3){(p-1)l_2}
\ar@{-}[d]^-{(p-1)l_1}
\\
\framebox{$\begin{array}{c}
      (12)
\end{array} : \begin{array}{l}
-(p-1)((p+1)l_1^2+l_2^2+(p+2)l_{1}l_2)=       \\
(p-1)(pl_2^2-l_{1}l_2) \\ \hline   [\Rightarrow l_1^2+l_2^2+l_{1}l_2 = 0] 
\end{array}  $} 
\ar@{-}[rd]^(.65){(p-1)l_1}
&
&
\framebox{$\begin{array}{c}
      (23)
\end{array} :  \begin{array}{l} (p-1)(pl_{1}l_2 - l_1^2)= \\ 
(p-1)(pl_1^2+(p+1)l_2^2+(2p+1)l_{1}l_2)
\\ \hline   [\Rightarrow l_1^2+l_2^2+l_{1}l_2 = 0]
\end{array}$ }
\ar@{-}[ld]
\ar@<-3pt>@{}[ld]_(.7){(p-1)l_2}
\ar@{-}[d]^-{-(p+1)l_1-pl_2}
\\
& 
\framebox{$1 :\begin{array}{l}
     (p-1)(p^2-1)l_{1}l_2^2  
\end{array}$}
&
\framebox{$\emptyset :  \begin{array}{l} -p(p^2-1)(l_1^{2}l_2+l_{1}l_2^2)=0 \\ \hline [\Rightarrow l_1^{2}l_2+l_{1}l_2^2 = 0]
\end{array}$
}
}
$
}
\end{figure}
Similarly to the unitary $\U(2,1)$ inert case, the inequalities gotten from the existence of partial Hasse invariants cutting out codimension 1 strata-closures are insufficient to prove that every tautological generically ordinary curve is strata-effective. In contrast with the inert case, the nefness of the Hodge line bundle isn't enough to additional information to conclude strata-effectivity. We haven't been able to determine whether (a) there is some tautological generically ordinary curve that is \textit{not} strata-effective or (b) every tautological generically ordinary curve is strata-effective but this cannot be detected by our methods.
\subsection{Guide to diagrams}
The boxes are indexed by elements of the Weyl group, which are displayed on the left-hand side. The cycle class of the corresponding stratum closure $[\overline{Y}_{w}] \in \Tsf^*(Y)$ is displayed on the right-hand side. A line from one box to another represents a partial Hasse invariant on the source stratum-closure which cuts out the target stratum-closure. It is labelled with the chern class of the line bundle of which the partial Hasse invariant is a section, expressed in $\Tsf^*(Y)$. To compute the cycle class of a stratum-closure one can choose an arrow into it and multiply the cycle class of the source stratum-closure by the chern class of the partial Hasse invariant. Where there are multiple paths to a stratum-closure this gives a relation in $\Tsf^*(Y)$ which is written under the line.
\subsection{Non-nefness of Hodge vector bundle}
\label{subsec-non-nef-A2}
One can use diagram \ref{fig-C2-min-taut}to prove that the Hodge vector bundle $\Omega^{\can}$ on the toroidal compactification of the Siegel-type Shimura variety for $g=2$ is not nef. This is in contrast to the situation in characteristic $0$.

\begin{proposition}
\label{prop-hodge-not-nef}
Assume that $\type(G)=\Csf_2$ and that $\mu^{\ad}$ is minuscule. Suppose that $X$ is projective and that $\zeta$ is smooth and surjective. Then the Hodge vector bundle $\Omega^{\can}$ is not nef on $X$.    
\end{proposition}
By~\ref{sec-intro-Hodge-type}, an immediate consequence is:
\begin{corollary} Let $S^{\Sigma}_{2,K'}$ be a smooth toroidal compactification of a Siegel threefold as in~\ref{sec-intro-Hodge-type}.  The canonical extension $\Omega^{\can}$ of the Hodge vector bundle on $S^{\Sigma}_{2,K'}$ is not nef.
\end{corollary}
\subsubsection{Proof of~\ref{prop-hodge-not-nef}}
    Here $G^{\ad} \cong \Sp(4)^{\ad}$ and the assumptions of~\ref{th-gen-ord-Sp4} hold. We use \ref{fig-C2-min-taut} to prove the statement. It suffices to give some curve $C \subset \PP(\Omega^{\can})$ such that $c_{1}(\mathcal{O}(1)).[C]$ has strictly negative degree. For $g=2$, $\PP(\Omega^{\can}) \cong Y_{2,K}^{\Sigma}$ is the flag space and $\mathcal{O}(1) \cong \mathcal{L}(0,-1)$. This implies that $c_{1}(\mathcal{O}(1)) = l_{2}$ in the notation of \ref{fig-C2-min-taut}. In fact it follows from \ref{fig-C2-min-taut}, that 
    \begin{equation}
        l_{2}.[\overline{Y}_{\sgn_{2}}] = -(p^2-1)l_{1}l_{2}^{3}
    \end{equation} has strictly negative degree and
    \begin{equation}
        l_{2}.[\overline{Y}_{(12)}] = (p-1)(p^2+1)l_{1}l_{2}^{3}
    \end{equation} has strictly positive degree. Hence, the line bundle $\mathcal{O}(1)$ is neither nef nor anti-nef. In particular, this means that $\Omega$ is not nef.
\qed

\section{Effective tautological classes III: Linear stratifications}
\label{sec-linear}
Let $(G,\mu)$ be a cocharacter datum. Let $d:=|\Phi \setminus \Phi_I|/2$.  For $0\leq i \leq d$, define $M_{i} := |\{w \in {}^{I}W \mid l(w)=i\}|$ to be the number of elements of $\iw$ of length $i$. 
\begin{definition}
\label{def-linear-strat}
The type of the Zip stratification of $\GZip^{\mu}$ is the sequence $(M_0,\ldots,M_d)$. The Zip stratification is linear if its type is $(M_0,\ldots,M_d) = (1,\ldots,1)$.
\end{definition}
\subsubsection{\texorpdfstring{$\textdbend$}{}} Let $\zeta: X \rightarrow \GZip^\mu$ be a Zip period map. The stratification of $X$ may still have several components of strata in each dimension even if the stratification is linear. For example, if $\zeta$ is not smooth then the strata $X_w$ need not be equi-dimensional.

\subsection{Linear Zip stratifications
}

\begin{proposition}
\label{prop-linear-classification}
The Zip stratification is linear if and only if  $(\type(G), \type(L))=(\Xsf_n, \Xsf_{n-1})$ for some $\Xsf \in \{\Asf, \Bsf, \Csf\}$ or $(\type(G), \type(L))=(\Gsf_2, \Asf_{1})$.   
\end{proposition} 
 \subsubsection{Proof}
The stratification is linear if and only if for every $l$, $0 \leq l \leq |\Phi^+ \setminus \Phi_I^+|$ there exists a unique $w \in \iw$ of length $l$. So this happens if and only if 
 \addtocounter{equation}{-1}
 \begin{subequations}
 \begin{equation}
 \label{eq-strat-linear}
 |\iw|=1+\frac{|\Phi \setminus \Phi_I|}{2}.    
 \end{equation}    
 \end{subequations}
 Two observations reduce the classification of the pairs $(\Delta, I)$ satisfying~\eqref{eq-strat-linear} to a small number of computations which are easily read off from Bourbaki's Planches \cite{bourbaki-lie-4-6}. 
 \subsubsection{\texorpdfstring{$L$}{} is Maximal}
 If the Zip stratification is linear, then $L$ is necessarily maximal: If $\alpha, \beta \in \Delta \setminus I$ are two distinct simple roots then the root reflections $s_{\alpha}$ and $s_{\beta}$ are two distinct length-one elements in $\iw$. Write $\Delta=I \sqcup \{\alpha\}$.
\subsubsection{Detecting non-linearity from a parabolic subgroup} Recall \cite[\S5.5]{Humphreys-coxeter-book} that for all $A \subset \Delta$, the Bruhat-Chevalley order (resp. length function) of the standard Coxeter parabolic subgroup $W_A \subset W$ is the restriction of the Bruhat-Chevalley order (resp. length function) of $W$ to $W_A$. It follows that, if $A \subset \Delta$ contains $\alpha$ and there are two elements of the same length in $\leftexp{I\cap A}{W_A}$, then the original Zip stratification is not linear. 
This reduces the problem to showing that~\eqref{eq-strat-linear} fails for the following pairs $(\type(G), \type(L))$: $(\Asf_3, \Asf_1 \times \Asf_1)$, $(\Csf_3, \Asf_2)$, $(\Bsf_3, \Asf_2)$, $(\Fsf_4, \Bsf_3)$, $(\Fsf_4, \Csf_3)$, $(\Esf_6, \Dsf_5)$, $(\Esf_7, \Esf_6)$ and $(\Esf_8, \Esf_7)$. In these cases, both sides of~\eqref{eq-strat-linear} are swiftly computed using the Planches \cite{bourbaki-lie-4-6} and one sees that the left-hand side is strictly larger.

 \qed
\subsubsection{Proof of \texorpdfstring{~\ref{th-linear-eff}}{}}
    Every graded piece of $\CHQ^*(\GZip^\mu)$ is $1$-dimensional since the Zip stratification is linear. Hence, every graded piece of $\Tsf^*(X)$ is $1$-dimensional and so any $\eta \in \Tsf^i(X)$ is of the form $\eta = a_{\eta}[\overline{X}_{w}]$ for $l(w) = d-i$ where $a_{\eta} \in \QQ$. It remains to check that $a_{\eta}$ is non-negative if $\eta$ is effective. 
By the classification~\ref{prop-linear-classification} $\mu$ is quasi-constant unless $\type(G)=\Gsf_2$ and then $\mu$ is orbitally $5$-close.
Thus the hypotheses imply that the Griffiths line bundle $\grif(\GZip^{\mu},r)$~\ref{sec-grif-GZip} admits a Hasse invariant on every stratum closure. Hence 
    \begin{equation*}
        a_{\eta}c_1(\zeta^*\grif(\GZip^{\mu},r))^{l(w)}.[\overline{X}_w] = c_1(\zeta^*\grif(\GZip^{\mu},r))^{l(w)}.\eta 
    \end{equation*}  is non-negative.
\qed

Given $X \rightarrow \GZip^{\mu}$ we have defined a notion of strata-effective classes on both $X$ and its flag space $Y/X$. In section \ref{sec-cone-conj} we linked strata-effectivity of classes on $Y$ with a certain cone conjecture. Clearly, if a class $\alpha \in \CHQ^*(Y)$ is strata-effective then $\pi_{*}(\alpha) \in \CHQ^*(X)$ is strata-effective. The following example shows that we can have effective tautological classes on $Y$ which are not strata-effective while all effective tautological classes on $X$ are strata-effective.
\begin{remark}
    There exist cocharacter data $(G,\mu)$ such that the stratification of $\GZip^{\mu}$ is linear but $\Cscr_{\pha} \neq \Cscr_{\GZip^\mu}$.
\end{remark}
An example of this occurs for Shimura varieties with group $\U(n,1)$ of type $\Asf_n$.
\begin{corollary}
    There exist Zip period maps $\zeta : X \rightarrow \GZip^\mu$ such that every effective class in $\Tsf^*(X)$ is strata-effective but there is some effective $\eta \in \Tsf^1(Y)$ which is not strata-effective.
\end{corollary}

\bibliographystyle{plain}
\bibliography{biblio_overleaf}
\end{document}